\documentclass{elsart}
\usepackage{graphicx}
\newcommand{\ba}{{\bf a}}

\newcommand{\bxi}{{\mbox{\boldmath $\xi$}}}
\newcommand{\bep}{{\mbox{\boldmath $\epsilon$}}}

\newcommand{\bug}{{\underline {\bf g}}}
\newcommand{\bue}{{\underline {\bf e}}}
\newcommand{\buq}{{\underline {\bf q}}}
\newcommand{\buv}{{\underline {\bf v}}}
\newcommand{\buw}{{\underline {\bf w}}}
\newcommand{\bU}{{\bf U}}
\newcommand{\bP}{{\bf P}}
\newcommand{\bQ}{{\bf Q}}
\newcommand{\bG}{{\bf G}}
\newcommand{\bR}{{\bf R}}
\newcommand{\bA}{{\bf A}}
\newcommand{\bC}{{\bf C}}
\newcommand{\bE}{{\bf E}}
\newcommand{\bV}{{\bf V}}

\newcommand{\R}{I \! \! R}

\newcommand{\C}{I \! \! \! \! {C}}

\newcommand{\uu}{{\underline u}}

\newcommand{\ua}{{\underline a}}

\newcommand{\umu}{{\underline \mu}}

\newcommand{\ux}{{\underline x}}

\newcommand{\ue}{{\underline e}}

\newcommand{\p}{{p^*}}

\newcommand{\RZ}{{\mathcal{R}(Z)}}

\newcommand{\RA}{{\mathcal{R}(\bA_k)}}
\newcommand{\bb}{\begin{eqnarray}}
\newcommand{\be}{\end{eqnarray}}

\newtheorem{theorem}{Theorem}
\newtheorem{proposition}{Proposition}
\newtheorem{lemma}{Lemma}
\newtheorem{corollary}{Corollary}

\begin{document}

\title{On the condensed density of the generalized eigenvalues of pencils
of  Gaussian random matrices and applications}
\author{P. Barone }
\address{ Istituto per le Applicazioni del Calcolo ''M. Picone'',
C.N.R.,\\
Via dei Taurini 19, 00185 Rome, Italy \\
e-mail: piero.barone@gmail.com, p.barone@iac.cnr.it  \\
fax: 39-6-4404306}

\maketitle

\section*{Abstract}

Pencils of  matrices whose elements have a joint noncentral
Gaussian distribution with nonidentical covariance are
considered. An approximation to the distribution of the
squared modulus of their determinant is computed which
allows to get a closed form approximation of the condensed
density of the generalized eigenvalues of the pencils.
Implications of this result for solving several moments
problems are discussed and some numerical examples are
provided.

{\it Key words:}
 random determinants,
 complex exponentials, complex moments problem,
 logarithmic potentials

\newpage

\section*{Introduction}

Let us define the random Hankel matrices \bb
\bU_0=\left[\begin{array}{llll}
\ba_0 & \ba_{1} &\dots &\ba_{p-1} \\
\ba_{1} & \ba_{2} &\dots &\ba_{p} \\
. & . &\dots &. \\
\ba_{p-1} & \ba_{p} &\dots &\ba_{n-2}
  \end{array}\right],\;\;
\bU_1=\left[\begin{array}{llll}
\ba_1 & \ba_{2} &\dots &\ba_{p} \\
\ba_{2} & \ba_{3} &\dots &\ba_{p+1} \\
. & . &\dots &. \\
\ba_{p} & \ba_{p+1} &\dots &\ba_{n-1}
  \end{array}\right]\label{pencil}\be
where $n=2p,\;\;\ba_k=s_k+\bep_k,\quad
k=0,1,2,\dots,n-1,\;\;\bep_k$ is a complex Gaussian, zero mean,
white noise,
 with variance $\sigma^2$ and $s_k\in\C$.
In the following all random quantities are
denoted by bold characters.
 Let us consider the generalized eigenvalues $\{\bxi_j,\;j=1,\dots,p\}$ of $(\bU_1,\bU_0)$ i.e. the
 roots of the polynomial
$\bP(z)=\det[(\bU_1-z\bU_0)]$, which form a set of exchangeable random variables. Their marginal density $h(z)$ also called condensed density
\cite{ham} or normalized one-point correlation function \cite{deift}, is the expected value of the
(random) normalized counting measure on the zeros of $\bP(z)$ i.e.
$$h(z)= \frac{1}{p}E\left[\sum_{j=1}^{p}\delta(z-\bxi_j)\right]$$ or,
equivalently, for all Borel sets $A\subset\C$ $$\int_A
h(z)dz=\frac{1}{p}\sum_{j=1}^p Prob(\bxi_j\in A).$$ It can be proved
that (see e.g. \cite{barja}) $h(z)=\frac{1}{4\pi}\Delta u(z) $ where
$\Delta$ denotes the Laplacian operator  with respect to $x,y$ if
$z=x+iy$ and $u(z)=\frac{1}{p}E\left\{\log(|\bP(z)|^2)\right\}$ is the
corresponding logarithmic potential.

The condensed density $h(z)$ plays an important role for solving
moment problems  such as the trigonometric,  the complex, the
Hausdorff ones. It was shown in $[2-10]$, \cite{maba2,maba1} that all
these problems can be reduced  to the complex exponentials
approximation problem (CEAP), which can be stated as follows. Let us
consider  a uniformly sampled signal made up of a linear combination
of complex exponentials
\begin{eqnarray}s_k=\sum_{j=1}^\p c_j\xi_j^k.\label{modale}\end{eqnarray}
where $c_j,\xi_j\in\C.$ Let us assume to know an even number $n=
2p,\;p\ge \p$ of noisy samples
$$\ba_k=s_k+\bep_k,\quad k=0,1,2,\dots,n-1$$ where $\bep_k$ is a complex Gaussian,
zero mean, white noise, with finite known variance $\sigma^2$.  We
want to estimate  $\p,c_j,\xi_j,\;j=1,\dots,\p$, which is a well
known ill-posed inverse problem. We notice that, in the noiseless
case and when $p=\p$, the parameters $\xi_j$ are the generalized
eigenvalues of the pencil $(U_1,U_0)$ where now $U_0$ and $U_1$ are
built as in (\ref{pencil}) but starting from $\{s_k\}.$

From its definition it is evident that the condensed density
provides information about the location in the complex plane of the
generalized eigenvalues of $(\bU_1,\bU_0)$ whose estimation is the most
difficult part of CEAP. Unfortunately its computation is very
difficult in general. In \cite{baram} a method to solve CEAP was
proposed based on an approximation of the condensed density. An
explicit expression of $h(z)$ proposed by Hammersley \cite{ham} when
the coefficients of $\bP(z)$ are jointly Gaussian distributed was
used. The second order statistics of these coefficients in the CEAP
case were estimated by computing many Pade' approximants of
different orders to the $Z$-transform of the data $\{\ba_k\}$. This
last step was essential to realize the averaging that appears in the
definition of $h(z)$, which is the key feature to make the condensed
density a useful tool for applications. In fact in the noiseless
case $h(z)$ is a sum of Dirac $\delta$ distributions centered on the
generalized eigenvalues while, when the signal is absent ($s_k=0\;\;
\forall k$), it was proved in \cite{barja} that if $z=r
e^{i\theta}$, the marginal condensed density $h^{(r)}(r)$ w.r. to
$r$  of the generalized eigenvalues is asymptotically in $n$ a Dirac
$\delta$ supported on the unit circle $\forall \sigma^2$. Moreover
for finite $n$ the  marginal condensed density w.r. to $\theta$ is
uniformly distributed on $[-\pi,\pi]$. Therefore if the
signal-to-noise ratio (e.g. $SNR=
\frac{1}{\sigma}\min_{h=1,\p}|c_h|$) is large enough $h(z)$ has
local maxima  in a neighbor of each $\xi_j,\;j=1,\dots,\p$ and this
fact can be exploited to get good estimates of $\xi_j$. However
usually we have only one realization of the discrete process
$\{\ba_k\}$, hence we cannot estimate $h(z)$ by averaging. In
\cite{barja2} a stochastic perturbation method is proposed to
overcome this problem. It is based on the computation of the
generalized eigenvalues of many pencils obtained by suitable
perturbations of the measured one. The computational burden is
therefore relevant. We then look for an approximation of $h(z)$
which can be well estimated by a single realization of $\{\ba_k\}$. It
turns out that the proposed approximation holds for pencils made up
by random matrices whose elements have a joint Gaussian
distribution. However the specific algebraic structure of CEAP,
which gives rise to pencils of Hankel random matrices, can be taken
into account to further reduce the computational burden. Moreover it
will be shown  that the noise contribution to $h(z)$
can be smoothed out to some extent simply acting on a parameter of
the approximant.

The paper is organized as follows. In Section 1 some
algebraic and statistical preliminaries are developed. In Section 2 the
closed form approximation of $h(z)$ is defined in the
general case.  In
Section 3 it is shown how to get a smooth estimate of
$h(z)$ from the data by exploiting its closed form
approximation in the Hankel case. In Section 4
computational issues are discussed in the Hankel case.
Finally in Section 5 some numerical examples are provided.

\section{Preliminaries}

Let us consider the $p\times p$ complex random pencil
$\bG(z)=\bG_1-z\bG_0,\;\;z\in\C$ where the elements of
$\Re{\bG_0},\Im{\bG_0},\Re{\bG_1},\Im{\bG_1}$ have
   a joint Gaussian
  distribution and $\Re$ and $\Im$ denotes the real and imaginary
  parts. Dropping  the dependence on
  $z$ for simplifying the notations, let us define
$$\bG=[\bug_1,\dots,\bug_p],\;\;\;\bug=vec(\bG)=[
\bug_1^T, \bug_2^T,\dots,\bug_{p}^T]^T
 $$
 Moreover $\forall z$, let us define $\check{\bug}_k
  =[
\Re{\bug_k}^T,\Im{\bug_k}^T]^T$
  and
 $$ \check{\bug}=[
\check{\bug}_1^T,\check{\bug}_2^T,\dots,\check{\bug}_{p}^T]^T.$$
   Then $\check{\bug}$ will have a multivariate
  Gaussian distribution with mean
  $\umu=E[\check{\bug}]\in\R^{2p^2}$
   and covariance $\Sigma\in\R^{2p^2\times 2p^2}$. We notice that
  no independence assumption neither between elements
  of $\bG_0$ and $\bG_1$
  nor between real and imaginary parts is made.
  Hence this is the most general hypothesis that can be done
  about the
  Gaussian distribution of the complex random vector $\bug$,
  (see \cite{vdbos} for a full discussion of this point).

Let us consider the $QR$ factorization of $\bG$ where
$\bQ^H\bQ=\bQ\bQ^H=I_p$ where $H$ denotes transposition plus
conjugation, $\bR$ is an upper triangular matrix and $I_p$ is
the identity matrix of order $p$. We then have
$$|\det(\bG)|^2=|\det(\bQ\bR)|^2=|\det(\bR)|^2=\prod_{k=1,p}|\bR_{kk}|^2.$$
We want to compute  the condensed density of the
generalized eigenvalues of the pencil $\bG(z)$ which
is given by \cite{ham,barja}: $$h(z)=\frac{1}{4\pi p}\Delta
E\left\{\log(|det[\bG(z)]|^2)\right\}= \frac{1}{4\pi p}\Delta
\sum_{k=1}^pE\left\{\log|\bR_{kk}(z)|^2\right\} .$$ We are therefore interested on the distribution of
$|\bR_{kk}|^2,k=1,\dots,p$ in order to compute
$E[\log|\bR_{kk}|^2]$.

To
perform the $QR$ factorization of the random matrix $\bG$ we
can use the Gram-Schmidt algorithm.  If  $\bQ=[\buq_1,\dots,\buq_{p}]$
it is given in Table \ref{tav1}.
%\clearpage
\linespread{1.}
\begin{table}[!h]
\begin{center}
\begin{tabular}{|l|}
\hline
For $k=1,\dots,p$ \\
\hspace{0.4cm}$\buw_k=\bug_k$ \\
\hspace{0.4cm}if $k>1$ then \\
\hspace{0.8cm}$\bR_{ik}=\buq_i^H\bug_k,\;i=1,\dots,k-1$\\
\hspace{0.8cm}$\buw_k=\buw_k-\sum_{i=1}^{k-1}\bR_{ik}\buq_i$\\
\hspace{0.4cm}end\\
\hspace{0.8cm}$\bR_{kk}=\sqrt{\buw_k^H\buw_k}$\\
\hspace{0.8cm}$\buq_k=\frac{\buw_k}{\bR_{kk}}$\\
end\\
\hline
\end{tabular}
\caption{The Gram-Schmidt algorithm}
\end{center}
 \label{tav1}
\end{table}
\linespread{1.5}

We notice that $|\bR_{kk}|=\bR_{kk}$ and  $$\bR_{kk}^2= \left\{\begin{array}{llll}
 \bug_k^H\bug_k, \mbox{ if $k=1$}\\
\bug_k^H\left(I_p-\sum_{i=1}^{k-1}\buq_i\buq_i^H \right)\bug_k, \mbox{ if $k>1$}
  \end{array}\right.$$
where  $\buq_i$ are functions of $\bug_j, \;j=1,\dots,i$. Therefore, denoting by $\tilde{\bug}_k=\{\check{\bug}_1,\dots,\check{\bug}_{k-1}\}$ we have that
$$\left\{\begin{array}{llll} \bR_{11}^2 \mbox{ is a quadratic form in Gaussian variables}\\
\bR_{kk}^2,\;k>1, \mbox{ conditioned on $\tilde{\bug}_k$, is a quadratic form in Gaussian variables }.
\end{array}\right.$$
Moreover
let us denote by $\ue_k$ the $k-$th column of $I_p$ and  let be $E_k=\ue_k\otimes I_{2p}$ then $\umu_k=E_k^T\umu,\;\;\Sigma_k=E_k^T\Sigma E_k$ are the mean vector and covariance matrix of $\check{\bug}_k$. Then we have

\begin{lemma}
For $k=1$ and for $k>1$, conditioned on $\tilde{\bug}_k$, $\bR_{kk}^2$ is distributed as $\sum_{r=1}^{n} \lambda^{(k)}_r \chi^2_{\nu_r}(\delta_r)$, $n=2p$, and $\chi^2_{\nu_r}(\delta_r)$ are independent, where $2(p-k+1)=\sum_{r=1}^n\nu_r\;$, $\lambda^{(k)}_r$  are the distinct eigenvalues  of $\Sigma_k^{1/2}{\mathcal{R}(\bA_k)}\Sigma_k^{1/2}\;\;$ with multiplicity $\nu_r$,  $\uu^{(k)}_i,\;i=1,\dots,n$ are the corresponding eigenvectors,
$\delta_r=\sum_{(r)}( \uu^{(k)}_i)^T\Sigma_k^{-1/2} \umu_k)^2,$ the summation being over all eigenvectors corresponding to eigenvalue $\lambda^{(k)}_r$, $$\bA_k=(I_p-\sum_{i=1}^{k-1}\buq_i\buq_i^H)$$ and
$${\mathcal{R}(\bA_k)}=\left[\begin{array}{rrrr}
\Re(\bA_k) &\;\; -\Im(\bA_k) \\
\Im(\bA_k) &\;\; \Re(\bA_k) \end{array}\right]$$
with $-\Im(\bA_k)=\Im(\bA_k)^T$ is the real isomorph of $\bA_k.$
\label{qf}
\end{lemma}

\noindent\underline{proof.} $\bA_k=\Re(\bA_k)+i\Im(\bA_k)$ is hermitian and idempotent because $\buq_i$ are orthonormal vectors. Therefore $\mbox{rank}(\bA_k)=p-k+1$ and $\mbox{rank}({\mathcal{R}(\bA_k)})=2(p-k+1)$ because the eigenvalues of ${\mathcal{R}(\bA_k)}$
are those of $\bA_k$ with multiplicity $2$. As $\bA_k$ depends only on $\buq_i,\;i=1,\dots,k-1$ which in turn depend only on $\bug_i,\;i=1,\dots,k-1$  it follows that, conditioned on $\tilde{\bug}_k$,  ${\mathcal{R}(\bA_k)}$ is a constant matrix. Moreover $\bug_k^H \bA_k \bug_k=\check{\bug}_k^T {\mathcal{R}(\bA_k)} \check{\bug}_k$ as can be easily checked. Let us define the variables
$\ux=(U^{(k)})^T\Sigma^{-1/2}\check{\bug}_k$. We have $\check{\bug}_k^T {\mathcal{R}(\bA_k)} \check{\bug}_k=\ux^T \Lambda^{(k)}\ux$ where $\Lambda^{(k)}$ is the diagonal matrix of eigenvalues of $\Sigma_k^{1/2}{\mathcal{R}(\bA_k)}\Sigma_k^{1/2}$. Only $m_k=2(p-k+1)$ eigenvalues are not zero and we can assume that they are the first $m_k$. Therefore $\bR_{kk}^2=\sum_{i=1}^{m_k}\lambda^{(k)}_i\ux_i^2$ is a quadratic form in Gaussian vectors of dimension $m_k$. The thesis follows e.g. by \cite[ch.29,sec.4]{jk}. $\Box$

\begin{corollary}
If $\Sigma=I_{2p^2}$ and $\umu=0$,  $\bR_{kk}^2$ is distributed as $\chi^2_{2(p-k+1)}$ \label{cor1}
\end{corollary}
\noindent\underline{proof.} As $\Sigma=I_{2p^2}$ the eigenvalues  of $\Sigma_k^{1/2}{\mathcal{R}(\bA_k)}\Sigma_k^{1/2}$ are those of ${\mathcal{R}(\bA_k)}$ which are $1$ with multiplicity $2(p-k+1)$ and $0$ with multiplicity $2(k-1)$. As $\umu=0$, $\delta_i=0$. Remembering that the $\chi^2_1(\delta_i)$ appearing in the previous Lemma are independent, the thesis follows
by the additivity property of $\chi^2$ distribution . $\Box$

\noindent\underline{Remark.} The Corollary follows also by  Bartlett's decomposition of a i.i.d. zero mean Gaussian random matrix \cite{bart}.

\section{Closed form approximation of $h(z)$}
Unfortunately we cannot use the easy result stated in the Corollary because in the case of interest the matrix $\bG(z)$ has a mean different from zero and a covariance structure depending on $z$. By  Lemma \ref{qf} we know that $\bR_{11}^2$ is distributed as a linear combination of non-central $\chi^2$ distributions. It is known that this distribution
admits an expansion $\mathcal{L} \left(\alpha,\beta,\tau\right)$ in series of generalized  Laguerre polynomials \cite[ch.29,sec.6.3]{jk}
 and the series is uniformly convergent in $\R^+$. More specifically let us denote the generalized Laguerre polynomial of order $m$  by
 $$L_m(x,\alpha)=\frac{x^{-(\alpha-1)}e^x}{m!}\frac{\partial^m}{\partial x^m}(x^{m+\alpha-1}e^{-x})=\sum_{h=0}^m c_{hm}x^{h}$$ where $$c_{hm}=\frac{(-1)^h\Gamma(\alpha+m)}{h!(m-h)!\Gamma(\alpha+h)}.$$
Then, following \cite{tz}, we have
\begin{lemma}
The density function of $\bR_{11}^2$ is given by
$$ f_1(y)=b_0\frac{y^{\alpha-1}e^{-y/\beta}}{\beta^\alpha\Gamma(\alpha)}+
\frac{y^{\alpha-1}e^{-y/\beta}}{\beta^\alpha\Gamma(\alpha)}\sum_{m=1}^\infty b_m L_m(y/\tau,\alpha)=\mathcal{L} \left(\alpha,\beta,\tau\right)$$
where $\alpha$ and $\beta$ are such that the first two moments of $\bR_{11}^2$ are the same of the first two moments of the gamma distribution representing the leading term of the expansion. Moreover $b_m$ are univocally determined by the moments and $\tau$ is a free parameter. If $\lambda_{max}$ denotes the maximum eigenvalue of $\Sigma_1$,  when $\tau^{-1}>2(\beta^{-1}-(2\lambda_{max})^{-1})$ the series  $\mathcal{L} \left(\alpha,\beta,\tau\right)$ is uniformly convergent
$\forall y\in \R^+.$ If $\beta>\lambda_{max}$ then $\tau=\beta$ makes the series to converge uniformly, $b_0=1$ and $b_m$ are determined by the first $m$ moments of $\bR_{11}^2$.
\label{distr1}
\end{lemma}
\noindent\underline{proof.} The proof follows by the results given in \cite{tz} for the distribution of quadratic forms in central normal variables which hold true also in the non-central case as can be easily checked. $\Box$

We can  compute $E[\log(\bR_{11}^2)]$ by
\begin{lemma}
\begin{eqnarray*}& &E[\log(\bR_{11}^2)]= b_0[\log\beta+\Psi(\alpha)]\\&+&\sum_{m=1}^\infty b_m\sum_{h=0}^m c_{hm}
\frac{\Gamma (\alpha+h)} {\Gamma (\alpha)}\left(\frac{\beta}{\tau}\right)^h \left[\log \beta+\Psi
   (\alpha+h)\right]
\end{eqnarray*}
\end{lemma}
\noindent\underline{proof.}
By Lemma \ref{distr1}  the series $\mathcal{L} \left(\alpha,\beta,\tau\right)$ converges uniformly. Therefore term-by-term integration can be performed and the result follows by noticing that, for $h=0,1,\dots$
$$\frac{1}{\beta^\alpha}\int_0^\infty\log(y)\left(\frac{y}{\tau}\right)^h y^{\alpha-1}e^{-y/\beta}dy=\Gamma (\alpha+h)\left(\frac{\beta}{\tau}\right)^h \left[\log \beta+\Psi(\alpha+h)\right].$$
 $\Box$

We have then obtained a closed form expression for $E[\log(\bR_{11}^2)]$ as a function of the moments of $\bR_{11}^2$. By noticing that  the same result holds true for the distribution of
$\bR_{kk}^2$ conditioned on $\tilde{\bug}_k$, we  show now how to get  an approximation of $E[\log(\bR_{kk}^2)],\;k>1$.
\begin{theorem}
The density function $f_k(y)$ of $\bR_{kk}^2$ can be expanded in a uniformly convergent series of Laguerre functions
$$ f_k(y)=b^{(k)}_0\frac{y^{\alpha_k-1}e^{-y/\beta_k}}{\beta_k^{\alpha_k}\Gamma(\alpha_k)}+
\frac{y^{\alpha_k-1}e^{-y/\beta_k}}{\beta_k^{\alpha_k}\Gamma(\alpha_k)}\sum_{m=1}^\infty b^{(k)}_m L_m(y/\tau_k,\alpha_k).$$
When the parameter $\tau_k$, that controls the uniform convergence of the series,  can be chosen equal to $\beta_k$, then $b^{(k)}_0=1$ and $b^{(k)}_m,\;m=0,\dots,N$ depends on the first $N+1$ moments of $\bR_{kk}^2$. Moreover
\begin{eqnarray*}& &E[\log(\bR_{kk}^2)]= b^{(k)}_0[\log\beta_k+\Psi(\alpha_k)]+\\&&
\sum_{m=1}^\infty b^{(k)}_m\sum_{h=0}^m c_{hm}
\frac{\Gamma (\alpha_k+h)} {\Gamma (\alpha_k)}\left(\frac{\beta_k}{\tau_k}\right)^h \left[\log \beta_k+\Psi
   (\alpha_k+h)\right].
\end{eqnarray*}
If the series is truncated after $N+1$ terms, the approximation error \begin{eqnarray*}\eta_N^{(k)}&=&\left|E[\log(\bR_{kk}^2)]-\left(b^{(k)}_0[\log\beta_k+\Psi(\alpha_k)]\right.\right.+\\&&
\left.\left.\sum_{m=1}^N b^{(k)}_m\sum_{h=0}^m c_{hm}
\frac{\Gamma (\alpha_k+h)} {\Gamma (\alpha_k)}\left(\frac{\beta_k}{\tau_k}\right)^h \left[\log \beta_k+\Psi(\alpha_k+h)\right]\right)\right|
\end{eqnarray*}
   is bounded by
\begin{eqnarray*} K_1\frac{ \epsilon^{N+1}}{\alpha_k^2\beta_k^{\alpha_k}\Gamma(\alpha_k)}\left( {_2F_2}\left(\alpha_k, \alpha_k; 1 + \alpha_k, 1 + \alpha_k;K_2\right)+
G^{3,0}_{0,2}\left(-K_2\left| \stackrel{ 1-\alpha_k , 1-\alpha_k }{\scriptstyle 0 ,  -\alpha_k ,  -\alpha_k} \right.\right)\right)\end{eqnarray*}
where $K_1>0,\; K_2>0$  and $0<\epsilon<1$ are constants,  $_2F_2$ is a generalized hypergeometric function and $G^{3,0}_{0,2}$ is a Meijer's G-function.\label{thk}
\end{theorem}
\noindent\underline{proof.} By Lemma \ref{qf}, conditioned on $\tilde{\bug}_k$, $\bR_{kk}^2$ is distributed as a linear combination of non-central $\chi^2$ distributions. Therefore the results obtained for $\bR_{11}^2$ are true and the conditional density of $\bR_{kk}^2$ given $\tilde{\bug}_k$, denoted by $f_k(y|\tilde{\bug}_k)$, exists as a function of $L_2[\R^+]$. But denoting by $h(\tilde{\bug}_k;\tilde{\umu}_k,\tilde{\Sigma}_k)$ the Gaussian density of $\tilde{\bug}_k$, the joint density $f_k(y,\tilde{\bug}_k)$ of $\bR_{kk}^2$ and $\tilde{\bug}_k$ is uniquely specified as $$ f_k(y,\tilde{\bug}_k)=f_k(y|\tilde{\bug}_k)h(\tilde{\bug}_k)$$  and the marginal w.r. to  $\bR_{kk}^2$ is $$f_k(y)=\int_{\R^{2p(k-1)}} f_k(y|\tilde{\bug}_k)h(\tilde{\bug}_k)d\tilde{\bug}_k.$$
As the Laguerre functions form a complete system of $L_2(\R^+)$
it must exist a convergent - in $L_2(\R^+)$ - Laguerre series representing $f_k(y)$ whose coefficient are univocally determined by the moments of $f_k(y)$  which are finite \cite[Th.4.1]{abate} and given by
$$\gamma_m=\int_{\R^{2p(k-1)}}\gamma_m(\tilde{\bug}_k) h(\tilde{\bug}_k)d\tilde{\bug}_k,\;\;m=1,2,\dots $$ where $\gamma_m(\tilde{\bug}_k)$ are the moments of $\bR_{kk}^2|\tilde{\bug}_k$. We show now that this Laguerre expansion is uniformly convergent.
By Lemma \ref{distr1} for all $\tilde{\bug}_k$ it exists a  constant $0<\tau(\tilde{\bug}_k)<\infty$ such that the Laguerre expansion $\mathcal{L} \left(\alpha(\tilde{\bug}_k),\beta(\tilde{\bug}_k),\tau(\tilde{\bug}_k)\right)$ of $f_k(y|\tilde{\bug}_k)$ is uniformly convergent in $\R^+.$ But then it exists $\tau_k>0$ such that  $\tau_k^{-1}=\sup_{\tilde{\bug}_k}\tau^{-1}(\tilde{\bug}_k)<\infty$ and $\mathcal{L} \left(\alpha(\tilde{\bug}_k),\beta(\tilde{\bug}_k),\tau_k)\right)$ is also
uniformly convergent in $\R^+\;\;\forall \tilde{\bug}_k.$ But then for \cite[Th.2.3e]{hen1} it follows that the Laguerre expansion $\mathcal{L} \left(\alpha_k,\beta_k,\tau_k\right)$ of $f_k(y)$ is uniformly convergent in $\R^+.$
We can then integrate term-by-term and we get
\begin{eqnarray*}& &E[\log(\bR_{kk}^2)]= b^{(k)}_0\log\beta_k+\Psi(\alpha_k)+\\&&
\sum_{m=1}^N b^{(k)}_m\sum_{h=0}^m c_{hm}
\frac{\Gamma (\alpha_k+h)} {\Gamma (\alpha_k)}\left(\frac{\beta_k}{\tau_k}\right)^h \left[\Psi
   (\alpha_k+h)+\log \beta_k\right]+\\&&
\int_0^\infty \log(y) e_N(y)dy
\end{eqnarray*}
where $$e_N(y)=\frac{y^{\alpha_k-1}e^{-y/\beta_k}}{\beta_k^{\alpha_k}\Gamma(\alpha_k)}\sum_{m=N+1}^\infty b^{(k)}_m L_m(y/\tau_k,\alpha_k).$$
In \cite[eq.(31)]{tz} the bound
$$|e_N(y)|\le K_1 \epsilon^{N+1}\frac{y^{\alpha_k-1}e^{K_2 y}}{\beta_k^{\alpha_k}\Gamma(\alpha_k)}
,\;\;0<\epsilon<1,\;K_2=-\beta_k^{-1}+\frac{R}{\tau_k(1+R)},\;\;\epsilon<R<1$$
is given where $K_1>0, K_2$ are constants. But then
\begin{eqnarray*}\left|\int_0^\infty \log(y) e_N(y)dy\right|\le\int_0^\infty\left|\log(y) e_N(y)\right|dy \le\\
 K_1\frac{ \epsilon^{N+1}}{\beta_k^{\alpha_k}\Gamma(\alpha_k)}\int_0^\infty|\log(y)|y^{\alpha_k-1}
e^{K_2 y}dy=\\K_1\frac{ \epsilon^{N+1}}{\beta_k^{\alpha_k}\Gamma(\alpha_k)}\cdot\left(
\int_1^\infty\log(y)y^{\alpha_k-1}e^{K_2 y}dy-\int_0^1\log(y)y^{\alpha_k-1}e^{K_2 y}dy\right)=\\
G^{3,0}_{0,2}\left(-K_2\left| \stackrel{ 1-\alpha_k , 1-\alpha_k }{\scriptstyle 0 ,  -\alpha_k ,  -\alpha_k} \right.\right)+\frac{1}{\alpha_k^2} {_2F_2}\left(\alpha_k, \alpha_k; 1 + \alpha_k, 1 + \alpha_k;K_2\right)
.\;\;\Box
\end{eqnarray*}

\noindent \underline{Remark.} The bound on the error given above is of little use in  practice because the computation of the constants $K_1,K_2$ is quite involved as they depend on all moments. However, by Corollary \ref{cor1} we know that when $\Sigma=I_{2p^2}$ and $\umu=0$ the expansion terminates after the first term and $b^{(k)}_0=1$. Moreover from Theorem \ref{thk} we know that it can happen that the coefficients  $b^{(k)}_k, k=0,\dots,N$ are determined by the first $N+1$ moments only. Therefore by continuity we can conjecture that the first term of the expansion provides most of the information in the general case and therefore the truncation error should be small. This conjecture is strongly supported by numerical evidence as discussed in Section 5.

 We can now  prove the main theorem
\begin{theorem}
If
$\bQ(z)\bR(z)$ is the $QR$ factorization of $\bG(z)$,
$$u(z)=\frac{1}{p}E\{\log(|\det[\bG_1-z\bG_0]|^2)\}=\frac{1}{p}\sum_{k=1}^pE\left\{\log|\bR_{kk}(z)|^2\right\}=$$
$$  \frac{1}{p}\sum_{k=1}^p \left\{b^{(k)}_0(z)[\log\beta_k(z)+\Psi(\alpha_k(z))]+\right.$$
$$\left.\sum_{m=1}^\infty b^{(k)}_m(z)\sum_{h=0}^m c_{hm}\frac{\Gamma (\alpha_k(z)+h)} {\Gamma (\alpha_k(z))}\left(\frac{\beta_k(z)}{\tau_k(z)}\right)^h \left[\log \beta_k(z)+\Psi(\alpha_k(z)+h)\right]\right\}.
$$
Moreover $$u(z)\approx \tilde{u}(z)=\frac{1}{p}\sum_{k=1}^p[\log\beta_k(z)+\Psi(\alpha_k(z))]$$
and $$|u(z)-\tilde{u}(z)|\le\frac{1}{p}\sum_{k=1}^p\eta^{(k)}_1(z).$$
\label{teo3}
\end{theorem}
 \noindent {\bf proof:  }
 By Theorem \ref{thk} we can approximate the density function  $f_k(y)$ of $\bR_{kk}^2$ by the first term divided by $b^{(k)}_0$  of its Laguerre expansion i.e. by $$\frac{y^{\alpha_k-1}e^{-y/\beta_k}}{\beta_k^{\alpha_k}\Gamma(\alpha_k)}.$$ This is a consistent approximation because this normalized first term is a $\Gamma$ density with parameters $\alpha_k,\beta_k$. But then the corresponding approximation of
 $E[\log(\bR_{kk}^2)]$ is $\log\beta_k+\Psi(\alpha_k)$  and then we
have
$$\tilde{u}(z)=\frac{1}{p}\sum_{k=1}^p[
\log[\beta_k(z)]+\Psi[\alpha_k(z)].$$ The other statements are obvious consequences of Theorem \ref{thk}. $\Box$

In order to compute the condensed density $h(z)$ we have to take the Laplacian of $u(z)$.
As differentiation can be a very unstable process, when we make use of the first order approximation of $u(z)$, we can expect that  even a small  approximation error in $u(z)$ can produce a large error in $h(z)$. However, in practice we have to approximate the Laplacian by finite differences by defining a square grid over the region of $\R^2$  which the unknown complex numbers $\xi_j$ are supposed to belong to. This provides an implicit regularization method if the grid size is properly chosen as a function of the approximation error of $u(z)$. We have
\begin{theorem}
If $\sup_{\C}\|u(z)-\tilde{u}(z)\|\le\varepsilon$ and if $z=x+iy$ and  the Laplacian operator is approximated by
\begin{eqnarray*}&&\hat{\Delta}u(x,y)= \\ &&\frac{1}{\delta^2}\left[u(x-\delta,y)+u(x+\delta,y)+u(x,y-\delta)+u(x,y+\delta)-4u(x,y) \right] \end{eqnarray*}
on a square grid with mesh size $\delta$ where $\delta(\varepsilon)=C\varepsilon^{1/4}, \;C$constant, then $$\|\hat{\Delta}\tilde{u}-\Delta u\|=O(\varepsilon^{1/4})$$ and this is the best possible approximation achievable.
\end{theorem}
\noindent\underline{proof.} By Taylor expansion of  $u(x\pm\delta,y\pm\delta)$ about $(x,y)$ we get $|\hat{\Delta}u(z)-\Delta u(z)|=O(\delta^2)$ and $$sup_{\C}|\hat{\Delta}u(z)-\Delta u(z)|=\|\hat{\Delta}u-\Delta u\|=O(\delta^2).$$ But $u(z)=\tilde{u}(z)+\eta(z)$ hence
$$\|\hat{\Delta}\tilde{u}-\Delta u\|=\|\hat{\Delta}u-\Delta u-\hat{\Delta}\eta\|\le O(\delta^2)+\frac{5\varepsilon}{\delta^2}.$$ For fixed $\varepsilon$ this error becomes unbounded as $\delta\rightarrow 0$. However by choosing
$\delta(\varepsilon)$ such that $\delta(\varepsilon)\rightarrow 0$ and $\frac{\varepsilon}{\delta(\varepsilon)}\rightarrow 0$ for $\varepsilon\rightarrow 0$ we get $$\|\hat{\Delta}\tilde{u}-\Delta u\|\rightarrow 0 \;\;\mbox{as  }\varepsilon\rightarrow 0$$Looking for a mesh size of the form $\delta(\varepsilon)=C\varepsilon^a$ such that the terms $O(\delta^2)$ and $\frac{5\varepsilon}{\delta^2}$ are balanced, we get
$$O(C^2\varepsilon^{2a})=O\left(\frac{5\varepsilon}{C^2\varepsilon^{2a}}\right)$$
which implies $a=\frac{1}{4}.$ In \cite{groe} it is proved that this bound is the best possible for all  approximation errors $\eta(z)$ such that $\|\eta\|\le \varepsilon$. $\Box$

As a final remark we notice that for computing $h(z)$ we
could start from the real isomorph \bb
{\mathcal{R}(\bG)}=\left[\begin{array}{rrrr}
\bV_R &\;\; -\bV_I \\
\bV_I &\;\; \bV_R \end{array}\right] \in\R^{n \times n}.
\nonumber \be of $\bG$ instead than from $\bG$. The following
proposition holds (\cite[Theorem 5.1]{ham}):
\begin{proposition}If
$\bG=\bV_R+i\bV_I,\;\;\bV_R,\bV_I\in\R^{p\times p}$,
 then
$$|\det(\bG)|^2=\det({\mathcal{R}(\bG)})$$
  \label{prop1}
\end{proposition}
 Let ${\mathcal{R}(\bG)}=\check{\bQ}\check{\bR}$ be the
  $QR$ factorization of ${\mathcal{R}(\bG)}$. Then have
$$|\det(\bG)|^2=det{\mathcal{R}(\bG)}=\prod_{k=1,n}\check{\bR}_{kk}$$
and \bb h(z)&=&\frac{1}{4\pi p}\Delta
E\left\{\log(|det[\bG(z)]|^2)\right\}= \frac{1}{4\pi p}\Delta
\sum_{k=1}^nE\left\{\log\check{\bR}_{kk}(z)\right\}
\nonumber\\&=& \frac{1}{4\pi n}\Delta
\sum_{k=1}^nE\left\{\log\check{\bR}^2_{kk}(z)\right\}.\nonumber
\be It will be shown in Section 4 however that this
expression of $h(z)$  is not convenient from the
computational point of view.

\section{Smooth estimate of the condensed density in the Hankel case}

We want to show now that we can exploit the closed form expression
of the condensed density to smooth out the noise contribution to
$h(z)$. This allows us to get a good estimate of $\p$ and
$\xi_j,j=1,\dots,\p,$ - which is the nonlinear  most difficult part of
CEAP - from a single  realization of the measured process $\{\ba_k\}$.
We first notice that by approximating the density of $\bR_{kk}^2$ by
a $\Gamma$ density with parameters $\alpha_k,\beta_k$, the mean and variance of $\bR_{kk}^2$
are approximated respectively by $\alpha_k\beta_k$ and $\alpha_k\beta_k^2$ and, if $b^{(k)}_0=1$, we have exactly
$$\gamma_1=\alpha_k\beta_k,\;\;\gamma_2=\alpha_k\beta_k^2+(\alpha_k\beta_k)^2. $$
However we know, by the proof of Theorem \ref{thk}, that
\bb\gamma_m=\int_{\R^{2p(k-1)}}\gamma_m(\tilde{\bug}_k) h(\tilde{\bug}_k)d\tilde{\bug}_k,\;\;m,1,2,\dots \label{integ}\be where $\gamma_m(\tilde{\bug}_k)$ are the moments of $\bR_{kk}^2|\tilde{\bug}_k$. The first two of them are given by (\cite{maprov})
$$\gamma_1(\tilde{\bug}_k)=tr[(\Sigma_k+2\umu_k\umu_k^T){\mathcal{R}(\bA_k)}]$$
$$\gamma_2(\tilde{\bug}_k)=2tr[(\Sigma_k+2\umu_k\umu_k^T){\mathcal{R}(\bA_k)}\Sigma_k{\mathcal{R}(\bA_k)}]+\gamma_1^2(\tilde{\bug}_k)$$
where $\umu_k=E[\check{\bug}_k]$ and $\Sigma_k=cov(\check{\bug}_k)$.
When $\bG=\bG_1-z\bG_0$ is an Hankel matrix and the elements of $\bG_0$ and $\bG_1$ are normally distributed with variance $\sigma^2$ as stated in the Introduction, it is easy to prove that the covariance matrix of $\bug_k$ does not depend on $k$ and it is a tridiagonal matrix $Z$ with $1+|z|^2$ on the main diagonal and $-z$ and $\overline{z}$ on the secondary ones. If $z=x+iy$ it turns out that $\forall k$ the covariance matrix  of $\check{\bug}_k$ is $\Sigma_k=\sigma^2\RZ$ where $\RZ$ is a $2\times 2$ block tridiagonal matrix given by
$$\RZ=\left[\begin{array}{cc} (-x,1+|z|^2,-x) & (y,0,-y) \\ (-y,0,y) & (-x,1+|z|^2,-x) \end{array}\right]$$ where $(a,b,c)$ denotes a tridiagonal matrix with $b$ on the main diagonal, $a$ and $c$ on the lower and upper diagonals respectively.
But then we have
$$\gamma_1(\tilde{\bug}_k)=\sigma^2 tr[\RZ\RA]+2tr[\umu_k\umu_k^T\RA]$$
\begin{eqnarray*}\gamma_2(\tilde{\bug}_k)&=&2\sigma^4tr[ \RZ\RA\RZ\RA]+\\
&&4\sigma^2tr[(\umu_k\umu_k^T)\RA \RZ\RA]+\gamma_1^2(\tilde{\bug}_k)\end{eqnarray*}
where the dependence on $\tilde{\bug}_k$ is only in $\RA$. By performing the integration in eq. (\ref{integ}) we get
$$\gamma_1=\sigma^2 c +d,\;\;
\gamma_2=\sigma^4 a+\sigma^2 b+\gamma_1^2$$
where
\begin{eqnarray*}a&=&2\int_{\R^{2p(k-1)}} tr[ \RZ\RA \RZ\RA]h(\tilde{\bug}_k)d\tilde{\bug}_k\\ b&=&4\int_{\R^{2p(k-1)}}tr[(\umu_k\umu_k^T)\RA \RZ\RA]h(\tilde{\bug}_k)d\tilde{\bug}_k\\
c&=&\int_{\R^{2p(k-1)}}
tr[\RZ\RA]h(\tilde{\bug}_k)d\tilde{\bug}_k\\ d&=&2\int_{\R^{2p(k-1)}}tr[\umu_k\umu_k^T\RA]h(\tilde{\bug}_k)d\tilde{\bug}_k. \end{eqnarray*}
But then we have
\begin{theorem}
If the density of  $\bR_{kk}^2$ is approximated by a $\Gamma$ density with parameters $\alpha_k,\beta_k $ such that the first two moments of $\bR_{kk}^2$ coincides with those of the approximant then
 $\beta_k$ is a  nondecreasing function of $\sigma^2$ and $\alpha_k$ is a nondecreasing function of $\frac{1}{\sigma^2}$ if $\frac{\|\umu_k\|_2^2}{\sigma^2}>\frac{E[tr(\RZ\RA)]}{2E[\tilde{\umu}_k^T\RA \tilde{\umu}_k]}$ where $\tilde{\umu}_k=\frac{\umu_k}{\|\umu_k\|_2}$ and $E$ denotes expectation w.r. to $h(\tilde{\bug}_k)$.
\end{theorem}
\noindent\underline{proof.}
$$\beta_k=\frac{\gamma_2-\gamma_1^2}{\gamma_1}=\sigma^2\left[\frac{\sigma^2 a+
b} {\sigma^2c+d}\right ],\;\;\;\;
\alpha_k=\frac{\gamma_1^2}{\gamma_2-\gamma_1^2}=\frac{\left(\sigma^2 c+d
\right)^2}{\sigma^4 a+\sigma^2b}$$ Deriving these
expressions respectively with respect to $\sigma^2$ and to  $\rho=\frac{d}{\sigma^2}$ where $d$ is assumed  fixed and $\sigma^2$ is variable, we get
$$\frac{\partial{\beta_k}}{\partial{\sigma^2}}=\frac{bd+a\sigma^2(2d+c\sigma^2)}
{(d+c\sigma^2)^2},\;\;\;\frac{\partial{\alpha_k}}{\partial{\rho}}=
\frac{2ad^2(\rho+c)+bd(\rho^2-c^2))}{(ad+b\rho)^2}.$$ But
$Z$ and $\bA_k$ are positive semidefinite matrices. Therefore
$\RZ$ and ${\mathcal{R}(\bA_k)}$ are also positive semidefinite because their eigenvalues are the same of those of $Z$ and $\bA_k$ with multiplicity $2$. Remembering that if $X,Y$ are positive semidefinite matrices $tr(X Y)^n\ge 0,\;n=1,2,\dots$, it follows that $a\ge0,b\ge0,c\ge0,d\ge0$ because the expectation of a nonnegative quantity is nonnegative. It follows that $\frac{\partial{\beta_k}}{\partial{\sigma^2}}\ge0.$
Moreover $\frac{\partial{\alpha_k}}{\partial{\rho}}\ge 0$ if
$$\rho^2-c^2=\frac{4}{\sigma^4}E[tr(\umu_k\umu_k^T\RA)]^2-c^2=
\frac{4(\umu_k^T\umu_k)^2}{\sigma^4}E[\tilde{\umu}_k^T\RA \tilde{\umu}_k]^2-c^2>0.$$
The thesis follows by noticing that $prob\{\tilde{\umu}_k^T\RA \tilde{\umu}_k>0\}=1$. In fact  $\RA$ is a random projector and the random quadratic form in the deterministic nonzero vector $\tilde{\umu}_k$ can  be zero only if $\tilde{\umu}_k$ is orthogonal to the random eigenvectors $\buv_{2k-1},\dots,\buv_{2p}$ of $\RA$ corresponding to nonzero eigenvalues. As this event has probability zero, $E[\tilde{\umu}_k^T\RA \tilde{\umu}_k]>0$.
$\Box$

The idea is then to use the parameters $\beta_k$ as
smoothing parameters and $\alpha_k$ as signal-related
parameters. By fixing
$\beta_k=\sigma^2\beta,\;\forall k$ and taking
$\alpha_k=\frac{\gamma_{1k}}{\sigma^2\beta}$ the
variance of $\bR^2_{kk}(z)$ is controlled by $\beta$ and
$h(z)$ can be estimated by
\bb\hat{h}(z)\propto\sum_{k=1}^p\hat{\Delta}
\left(\Psi\left[\frac{\hat{\gamma}_{1k}(z)}{\sigma^2\beta}\right]\right)\label{est}
\be where $\hat{\Delta}$ is the discrete Laplacian and
$\hat{\gamma}_{1k}(z)$ is an estimate of $\gamma_{1k}(z)$. In the
following we assume that the value $\hat{R}_{kk}^2(z)$ -
obtained by the $QR$ factorization of a realization of
$\bG(z)$ corresponding to a given set of observations
$\{a_k\}$ - is an estimate of the mode of $\bR_{kk}^2(z)$ and
therefore
$$\hat{R}_{kk}^2(z)=\beta_k(\alpha_k-1)$$ (see e.g
\cite{jk}[ch.17]). Then we get
$$\frac{\hat{\gamma}_{1k}(z)}{\sigma^2\beta}=\left(
\frac{\hat{R}_{kk}^2(z)}{\sigma^2\beta}+1\right).$$ From a
qualitative point of view,  increasing $\beta$ has the effect to
make larger the support of all modes of $h(z)$ and to lower their
value because $h(z)$ is a probability density. Hence the
noise-related modes are likely to be smoothed out by a sufficiently
large $\beta$. However a value of $\beta$ too large can result in a
low resolution spectral estimate.

\section{Computational issues in the Hankel case}

To estimate $h(z)$ on a lattice we must compute the $QR$
factorization of $\bG(z)$ for all values $z$ in the lattice.
This requires $O(m^2p^3)$  flops  if the lattice is square
of size $m$. However we notice that
$\bG(z)=\bU_1-z\bU_0=\bU(\bE_1-z\bE_0)$ where \bb
\bU=\left[\begin{array}{llll}
\ba_0 & \ba_{1} &\dots &\ba_{p} \\
\ba_{1} & \ba_{2} &\dots &\ba_{p+1} \\
. & . &\dots &. \\
\ba_{p-1} & \ba_{p} &\dots &\ba_{2p-1}
  \end{array}\right]\in \C^{p\times p+1}\be does not depend on $z$
  and
$$\bE_0=[\bue_1\dots \bue_{p}]\in \C^{p+1 \times p},\;\;
\bE_1=[\bue_2\dots \bue_{p+1}]\in \C^{p+1 \times p}$$ and
$\bue_k$ is the $k-$th column of the identity matrix of
order $p+1$. If $\bU=\bQ\bR$ is the $QR$ factorization of $\bU$
where $\bQ\in \C^{p\times p}$ is unitary and $\bR\in
\C^{p\times p+1}$ is upper trapezoidal, then the $QR$
factorization of $\bG(z)$ can be obtained simply by reducing
to upper triangular form by unitary transformations the
Hessemberg matrix
$$\bC(z)=\bR(\bE_1-z\bE_0)\in \C^{p\times p}.$$ This is the only task that
must be performed for each $z$. By using Givens rotations
this can be performed in $O(p)$ flops. The total cost of
the $QR$ factorization of $\bG(z)$ in the lattice reduces
then to $O(m^2p+p^3)$ flops.

\noindent Finally we notice that if we start from
${\mathcal{R}(\bG)}=\check{\bQ}\check{\bR},\;\;$$\bC(z)$ is  a
$2\times 2$ block matrix with Hessemberg diagonal blocks
and triangular off-diagonal ones. Therefore it can not be
transformed to triangular form in $O(2p)$ flops.

\section{Numerical results}

In this section some experimental evidence of the claims made in the
previous sections is given.

\begin{figure}
\centering{
\includegraphics[bb=  100  -100   1200   1000, height=15cm,width=15cm]{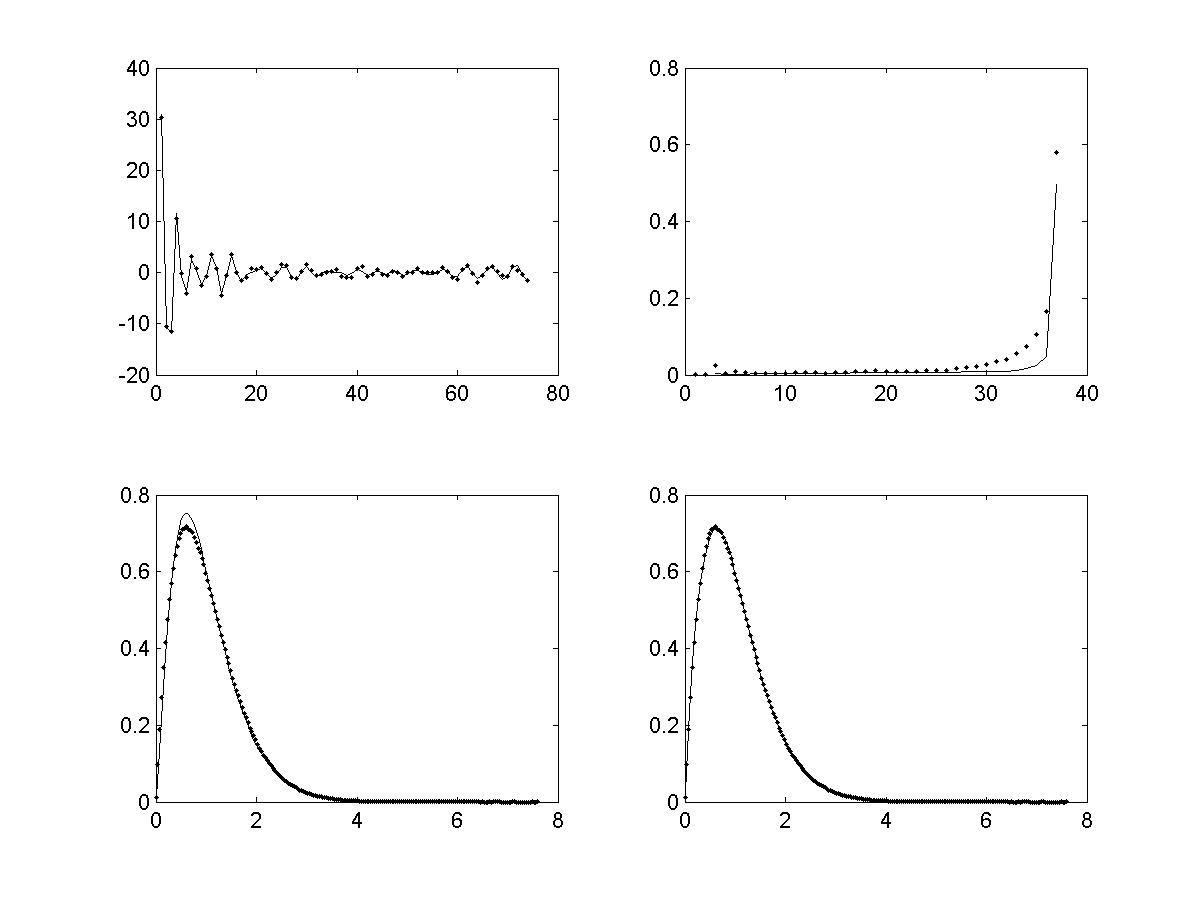}}
 \caption{Top left: real part of the signal (solid) and data (dotted) with $\sigma=0.5$;
 top right:   $L_2$ norm of the difference between the empirical density  of  $\bR_{kk}^2,k=1,\dots,30$  computed by MonteCarlo simulation with $4\cdot 10^6$ samples and its approximation obtained by truncating the series expansion of the density after the first term (dotted) and after the first $10$ terms (solid);
  bottom left:  density of $\bR_{kk}^2,\;k=36,$ approximated by the first term of its series expansion (solid), empirical density (dotted);
  bottom right: density of $\bR_{kk}^2,\;k=36,$ approximated by the first $10$ terms of its series expansion (solid), empirical density (dotted).
} \label{fig0}
\end{figure}

To appreciate the goodness of the approximation to the density of $\bR_{kk}^2$ provided by the truncated Laguerre expansion, $N=4\cdot 10^6$ independent
  realizations $a_k^{(r)},k=1,\dots,n,\;\;r=1,\dots,N$ of the r.v.  $\ba_k$ were
  generated from the complex exponentials model with $\p=5$ components given by
$$\underline{\xi}=\left[ e^{-0.1-i 2\pi  0.3},e^{-0.05-i 2\pi
0.28},e^{-0.0001+i 2\pi 0.2},e^{-0.0001+i 2\pi  0.21},e^{-0.3-i 2\pi
0.35}\right]$$ $$ \underline{c}=\left[ 6,3,1,1,20\right
],\;\;n=74,\;p=37,\;\sigma=0.5.$$ The matrices $U_0^{(r)},U_1^{(r)}$ based on $a_k^{(r)}$ were computed. The
 matrix $U_1^{(r)}-zU_0^{(r)}$ with $z=\cos(1)+i0.8$ was formed, its $QR$ decomposition  and the first $10$ empirical moments $\hat{\gamma}_j$ were computed. Estimates of the first $10$ coefficients of the Laguerre expansion were then computed by (\cite{balak})
 $$\hat{\alpha}_k= \frac{\hat{\gamma}_1^2}{\hat{\gamma}_2-\hat{\gamma}_1^2} ,\;\; \hat{\beta}_k=\frac{\hat{\gamma}_2-\hat{\gamma}_1^2}{\hat{\gamma}_1}  $$
 $$\hat{b}^{(k)}_h= (-1)^h\Gamma(\hat{\alpha}_k)\sum_{j=0}^h(-1)^j{h \choose j}\frac{\hat{\gamma}_{h-j}}{\Gamma(\hat{\alpha}_{k+h-j})},\;\;\hat{\gamma}_{0}=1,\;\;h=1,\dots,10$$
  The one term and ten terms approximations of the density were then computed and compared with the empirical density of $\bR_{kk}^2$ for $k=1,\dots,p$. The results are given in fig.\ref{fig0}. In the
top left part the real part of the signal  and of the data are plotted. In the
 top right part the   $L_2$ norm of the difference between the empirical density  of  $\bR_{kk}^2,k=1,\dots,p$  computed by MonteCarlo simulation and its approximation obtained by truncating the series expansion of the density after the first term and after the first $10$ terms is given.
 In the bottom left part the  density of $\bR_{kk}^2,\;k=36,$ approximated by the first term of its series expansion and the empirical density are plotted.
  In the bottom right part the density of $\bR_{kk}^2,\;k=36,$ approximated by the first $10$ terms of its series expansion and the empirical density are plotted.
  It can be noticed that the first order approximation is quite good even if it become worse for large $k$. The choice $\sigma=0.5$ is justified by the fact that this value is in the range of values used in the examples below. However the same kind of conclusions can be drawn for every SNR .

\begin{figure}
\centering{
\includegraphics[bb=  10  130   1200   1300,height=10cm,width=15cm]{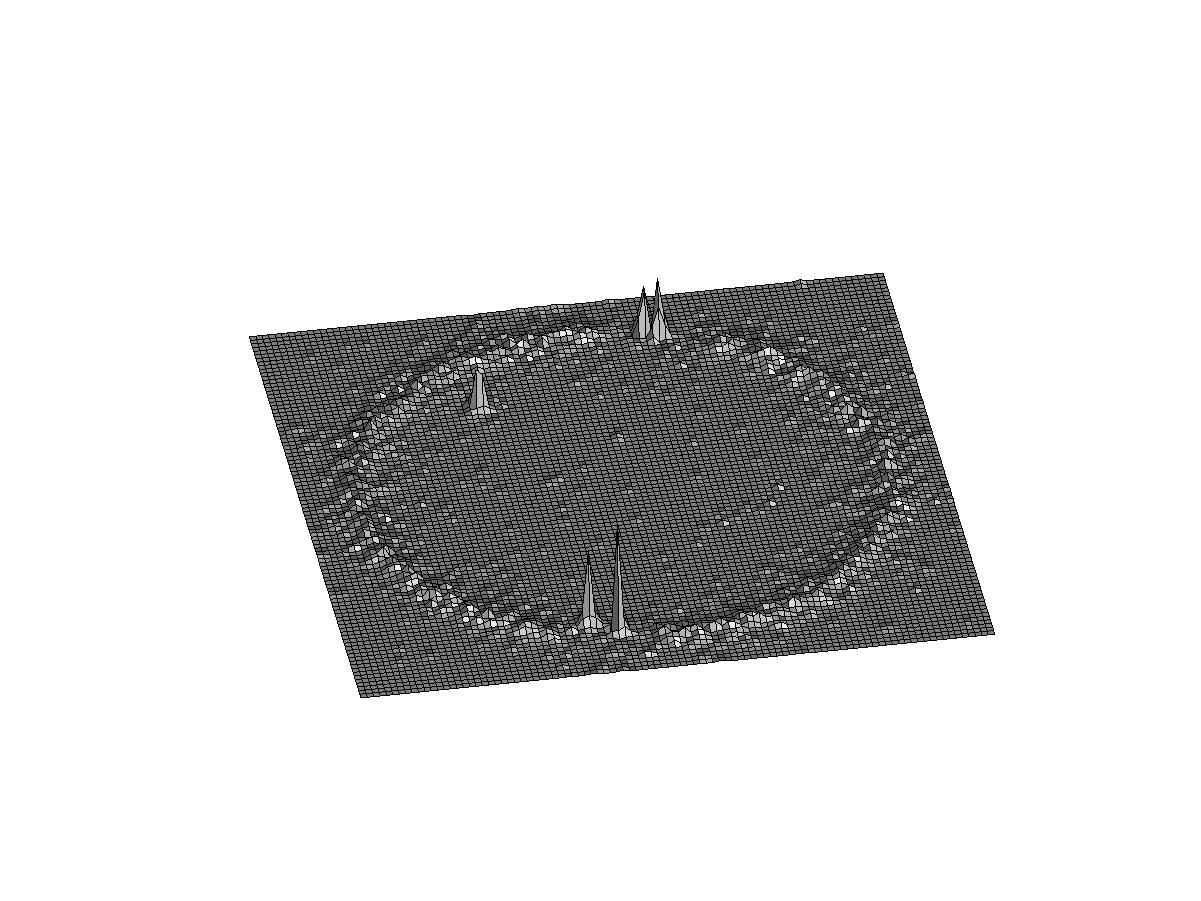}}
\includegraphics[bb=  10  130   1200   1300,height=10cm,width=15cm]{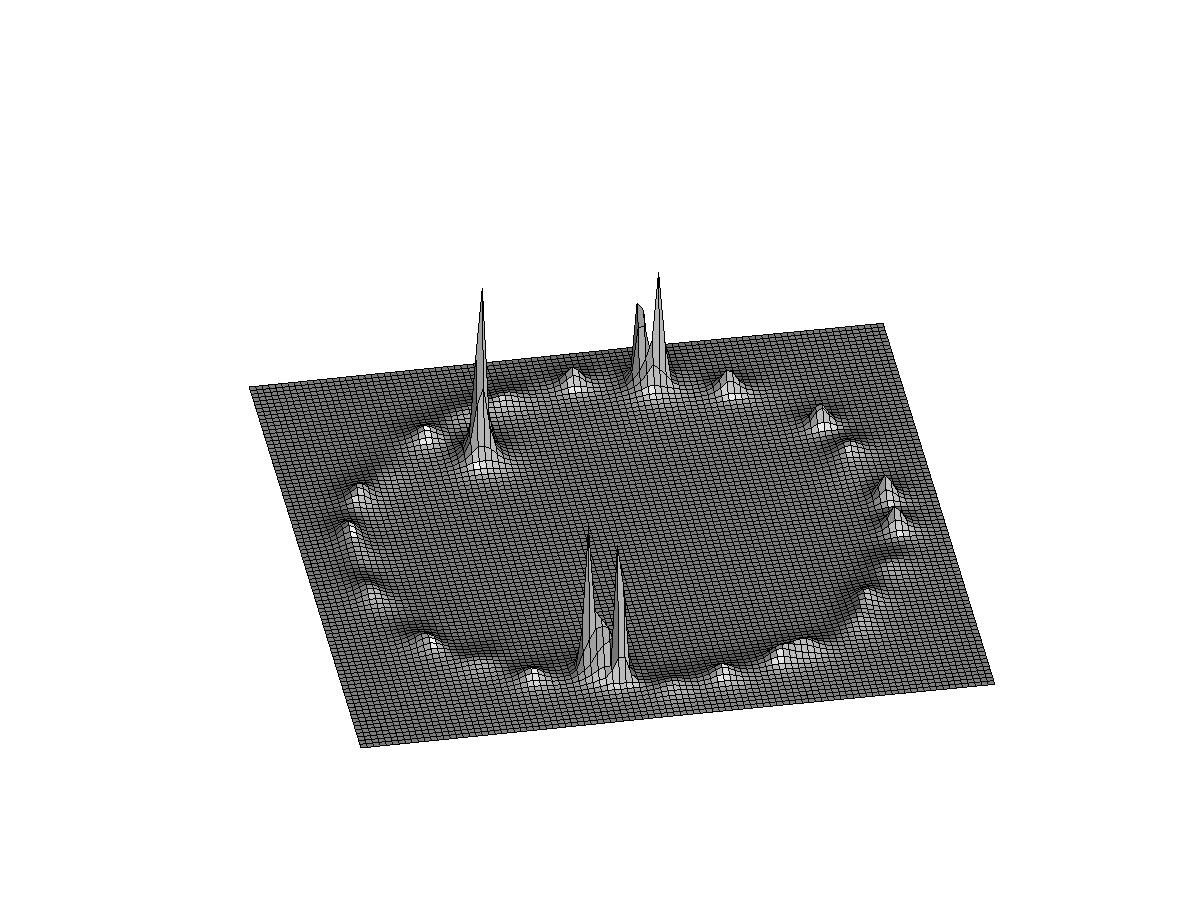}
 \caption{Top: Monte Carlo estimate of the condensed density when $\sigma=0.2$; bottom:
estimate of the condensed density by the closed form approximation
with $\beta=14.8$.} \label{fig11}
\end{figure}

\begin{figure}
\centering{
\includegraphics[bb= 10  130   1200   1300,height=10cm,width=15cm]{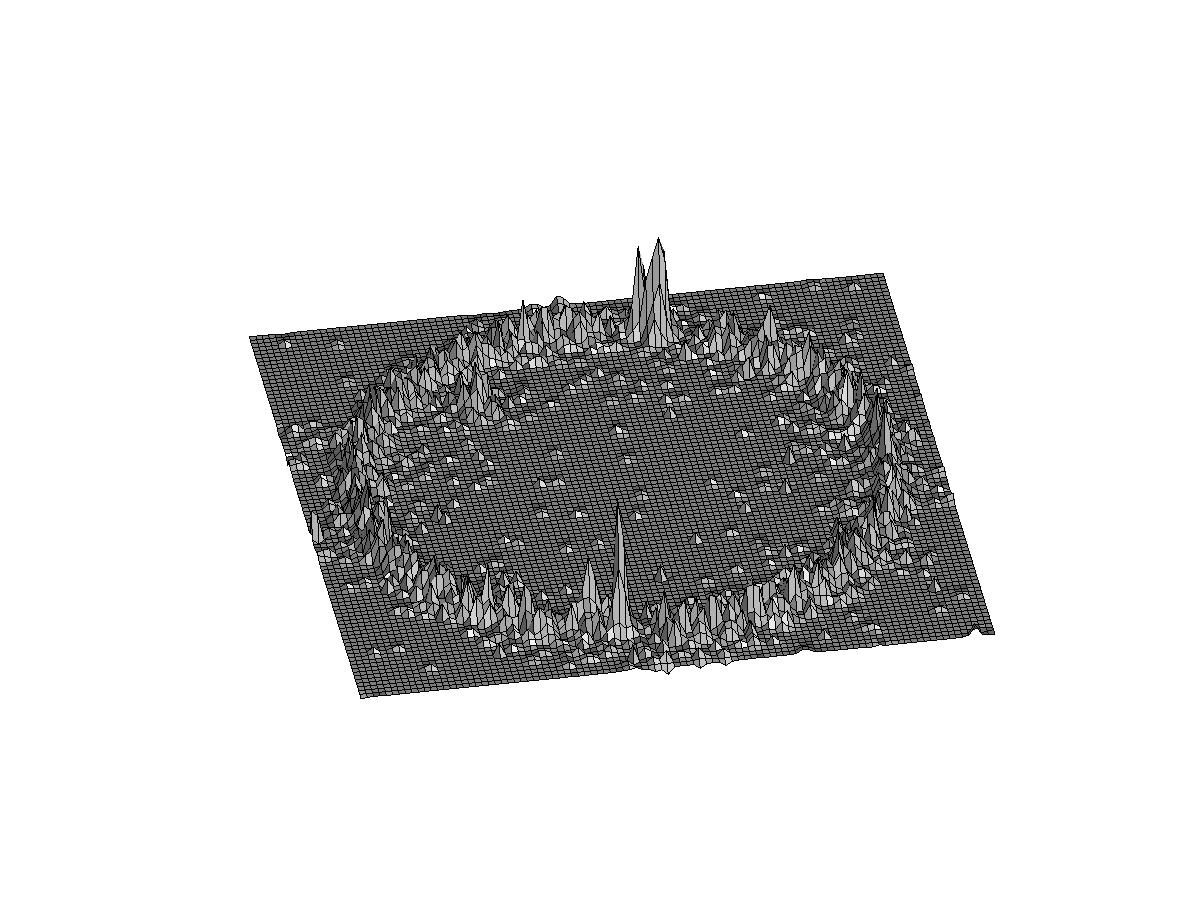}}
\includegraphics[bb= 10  130   1200   1300,height=10cm,width=15cm]{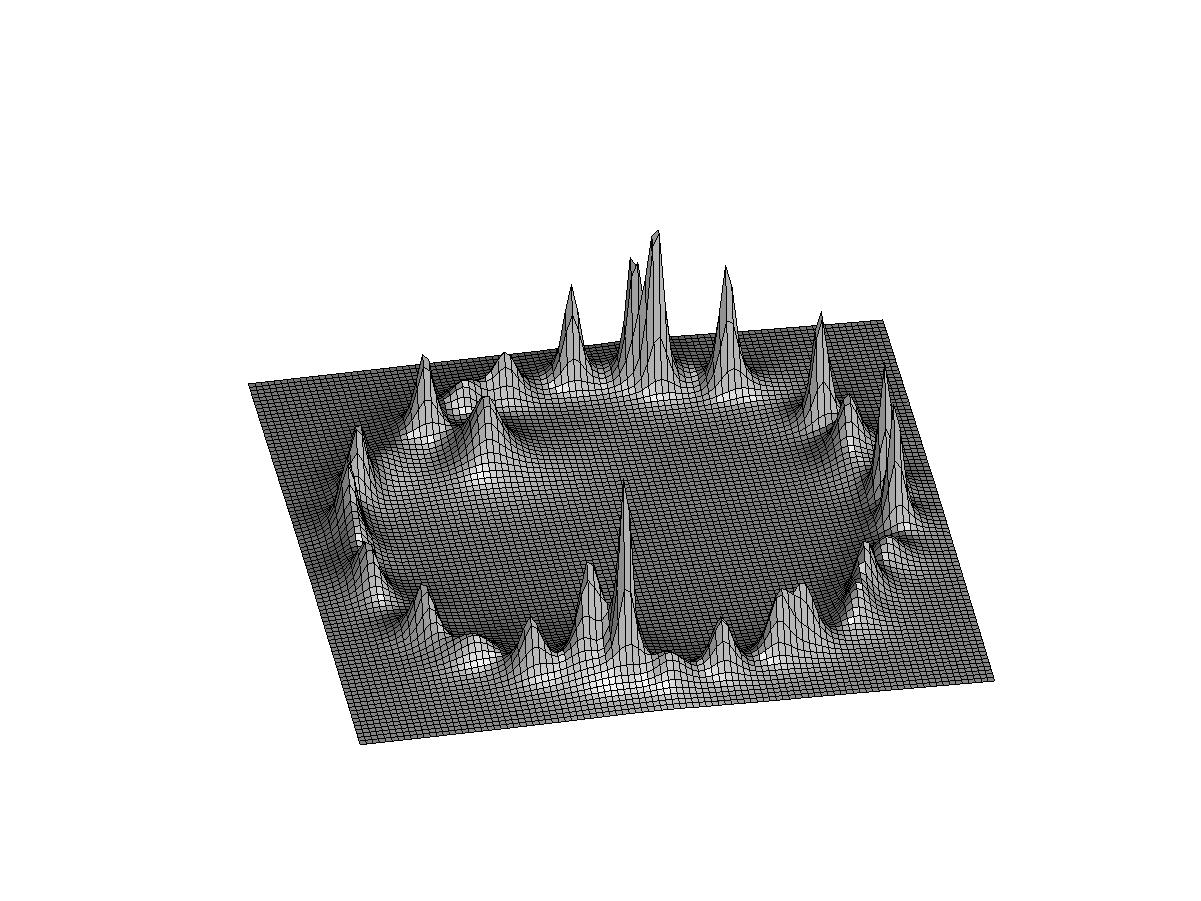}
 \caption{Top:  Monte Carlo
estimate of the condensed density when $\sigma=0.8$; bottom:
estimate of the condensed density by the closed form approximation
with $\beta=237$.} \label{fig12}
\end{figure}

To appreciate the advantage of the closed form estimate $\hat{h}(z)$
with respect to an estimate of the condensed density obtained by
MonteCarlo simulation an
 experiment was performed. $N=100$ independent
  realizations of the  r.v. generated above were considered. We notice that  the frequencies of the $3^{rd}$
and $4^{th}$ components are closer than the Nyquist frequency
($0.21-0.20=0.01<1/n=0.0135$). Hence a super-resolution problem is
involved in this case. Two values of the noise s.d. $\sigma$ were
used
$$\sigma=0.2,\;0.8.$$ An estimate of $h(z)$ was computed on a
square lattice of dimension $m=100$ by
$$\hat{h}(z)\propto
\sum_{r=1}^N\sum_{k=1}^p\hat{\Delta}
\left\{\Psi\left[\left(\frac{R_{kk}^{(r)}(z)^2}{\sigma^2\beta}+1\right)\right]\right\}
$$ where $R^{(r)}(z)$
is obtained by the QR factorization of the matrix $U_1^{(r)}-zU_0^{(r)}.$ In the top part of
fig.\ref{fig11} the estimate of $h(z)$ obtained by Monte
Carlo simulation   is plotted. In the bottom part the
smoothed estimates $\hat{h}(z)$ for $\sigma=0.2$ and
$\beta=5n$ based on a single realization was plotted. In
fig.\ref{fig12} the results obtained with $\sigma=0.8$ and
$\beta=5n$ are shown.  We notice that by the proposed
method we get an improved qualitative information with
respect to that obtained by replicated measures. This is an
important feature for applications where usually only one
data set is measured.  We also notice that when
$\sigma=0.2$ the probability to find a root of $P(z)$ in a
neighbor of $\xi_j$ is larger than the probability to find
it elsewhere. This is no longer true when $\sigma=0.8$ even
if the signal-related complex exponentials are well
separated. In the following we will say that the complex
exponential model is identifiable if this last case occurs
and it is strongly identifiable if the first case occurs.
Therefore if the model is identifiable
 the signal-related complex
exponentials are well separated  but the relative importance of some
of them - measured by the value of the local maxima of $h(z)$ - is
not larger than the relative importance of some noise-related
complex exponentials. Therefore in this case we need some a-priori
information about the location of the $\xi_j$ in order to separate
signal-related components from the noise-related ones.

We want now to show by means of a small simulation study the quality
of the estimates of the parameters $\p,\;\underline{\xi}$ and $
\underline{c}$ which can be obtained from $\hat{h}(z)$. To this aim
the following estimation procedure was used:
\begin{itemize}
\item the local maxima of $\hat{h}(z)$ are computed and sorted in
decreasing magnitude
\item a  clustering method  is used to group the
local maxima into two groups. If the model is strongly identifiable
  the signal-related maxima are larger than the
noise-related ones, therefore a simple thresholding is enough to
separate the two groups. A good threshold is the one that produces
an estimate of $s_k$ which best fits the data $a_k$ in $L_2$ norm as
the noise is assumed to be Gaussian
\item the cardinality $\hat{p}$ of the class with largest average value is an
estimate of $\p$
\item the local maxima $\hat{\xi}_j,j=1,\dots,\hat{p}$ of the class with
largest average value are estimates of $\xi_j,j=1,\dots,\p$. Of
course if $\hat{p}\ne \p$ some $\xi_j$ are not estimated or
viceversa some spurious complex exponentials are found
\item $
\underline{c}$ is estimated by solving the linear least squares
problem
$$\hat{\underline{c}}=argmin_{\ux}\|V\ux -\ua\|_2^2,\;\;\ua=[a_0,\dots,a_{n-1}]^T$$ where
$V\in\C^{n\times\hat{p}}$ is the Vandermonde matrix based on
$\hat{\xi}_j,j=1,\dots,\hat{p}$
\end{itemize}
The bias, variance and mean squared error (MSE) of each parameter
separately were estimated. $N=500$ independent data sets $\ua^{(r)}$
of length $n$ were generated by using the model parameters given
above and $\sigma=0.2$. For $r=1,\dots,N$  the condensed density
estimate $\hat{h}^{(r)}(z)$ was computed  on a square lattice of
dimension $m=100$. The estimation procedure is then applied to each
of the $\hat{h}^{(r)}(z),r=1,\dots,N$ and the corresponding
estimates
$\hat{\xi}_j^{(r)},\hat{c}_j^{(r)},j=1,\dots,\hat{p}^{(r)}$ of the
unknown parameters were obtained. If the estimate $\hat{p}^{(r)}$
was less than  the true value $\p$,  the corresponding data set
$\ua^{(r)}$ was discarded.

In Table \ref{tb2} the bias, variance and MSE of each parameter
including $\p$ is reported. They were computed by choosing among the
$\hat{\xi}_j^{(r)},j=1,\dots,\hat{p}^{(r)}\ge \p$ the one closest to
each $\xi_k,k=1,\dots,\p$ and the corresponding $\hat{c}_j^{(r)}$.
If more than one $\xi_k$ is estimated by the same
$\hat{\xi}_j^{(r)}$ the $r-$th data set $\ua^{(r)}$ was discarded.
In the case considered all the data sets were accepted.

As a second example the reconstruction of a piecewise constant
function from noisy Fourier coefficients is considered. The problem
is stated as follows. Given a real interval $[-\pi,\pi]$ and $N+1$
numbers $-\pi \leq l_1 < l_2 \dots < l_{N+1} \leq \pi$, let
$\mathcal{ F}$ be the class of functions defined as
$$F(t) = \sum_{j=1}^Nw_j\chi_j(t)\ ,$$
where
$$\chi_j(t) = \left\{ \begin{array}{ll}
 1  & \mbox{if $t\in [l_j,l_{j+1}]$} \\
 0 & \mbox{otherwise}
 \end{array}
 \right.\ ,$$
and the $w_j$ are real weights. The problem consists in
reconstructing a function $F(t)\in \mathcal{ F}$ from a finite
number of its noisy Fourier coefficients
$$\ba_k=\frac{1}{2}\int_{-\pi}^{\pi}F(t)e^{itk}dt+\bep_k =  s_k+\bep_k\ ,
\quad k=0,\dots,n-1\ ,$$ where ${\bep_k}$ is a complex Gaussian,
zero mean, white noise,
 with variance $\sigma^2$. We are looking for a solution which is not
affected by Gibbs artifact and can cope, stably, with the noise. The
basic observation is the following. The unperturbed moments $s_k$
are given by
$$s_k=\frac{1}{2}\int_{-\pi}^{\pi}F(t)e^{itk}dt=\sum_{j=1}^Nw_j
\frac{\sin(\beta_jk)}{k}\exp(i\lambda_jk),$$ where
\begin{eqnarray*}
\beta_j=\frac{l_{j+1}-l_j}{2},\;\; \lambda_j=\frac{l_{j+1}+l_j}{2}
.\end{eqnarray*} Then consider the  $Z$-transform of the sequence
$\{s_k\}$
$$s(z)=
\sum_{j=1}^Nw_j\left(\beta_j+\frac{1}{2i} \ln
\frac{z-e^{il_j}}{z-e^{il_{j+1}}}\right)$$ which converges
if $|z|>1$ and is defined by analytic continuation if
$|z|\le 1$. We notice that $s(z)$ has a branch point at
$\xi_j=e^{il_j},j=1,\dots,N+1$ where $l_j$ are the

discontinuity points of $F(t)$. It was proved in
\cite{maba2,maba1} that the $c_j$ are strong attractors of
the poles of the Pade' approximants $[q,r]_f(z)$ to the
noisy $Z$-transform
$$f(z)=\sum_{k=0}^\infty a_k z^{-k}$$ when $q,r\rightarrow\infty$ and
$q/r\rightarrow 1$. It is easy to show that the poles of
$[q,r]_f(z)$ are the generalized eigenvalues of the pencil
$(U_1,U_0)$ built from the data $a_k,k=0,\dots,{n-1}$ whose
condensed density is $h(z)$. Therefore, as shown in
\cite{maba2,maba1}  the local maxima of $h(z)$ are
concentrated along a set of arcs which ends in the branch
points $\xi_j$ and on a set of arcs close to the unit
circle. As the branch points are strong attractors for the
Pade' poles, the probability to find a pole in a neighbor
of a branch point is larger than elsewhere, therefore it
can be expected that the branch points correspond to the
largest local maxima of $h(z)$, as far as the SNR is
sufficiently large. In order to compute estimates
$\hat{l}_j$ of $l_j$, it is sufficient to compute the
arguments of the main local maxima of $\hat{h}(z)$. The
$w_j$ are then estimated by taking the median  in each
interval $[\hat{l}_j,\hat{l}_{j+1}]$ of the rough estimate
of $F(t)$ obtained by taking the discrete Fourier transform
of $a_k,\;k=0,\dots,n-1$. The median is in fact robust with
respect to errors affecting $\hat{l}_j$.

\begin{figure}[h]
\begin{center}
\centerline{
\includegraphics[bb= 10 130   1200   1300,height=7.5cm,width=8cm]{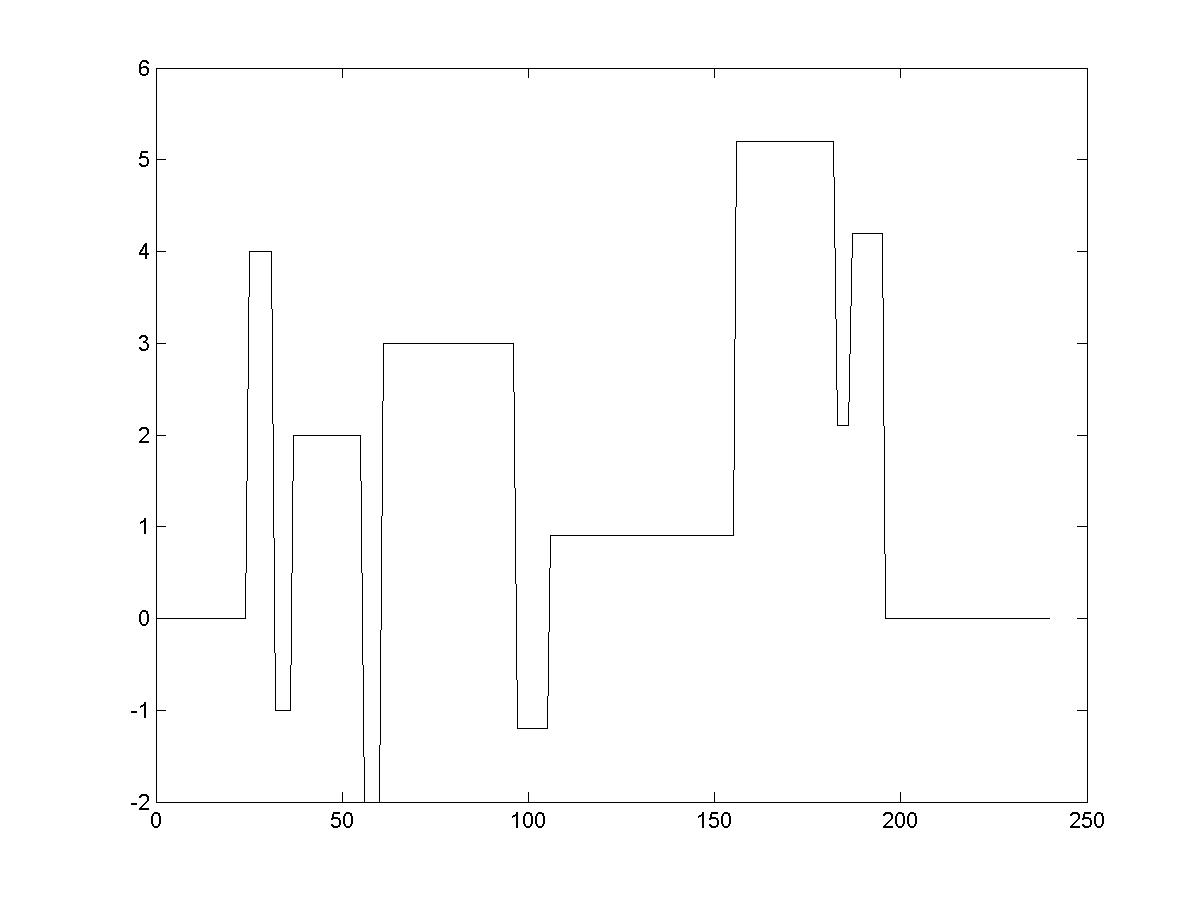}
\hspace{.5cm}
\includegraphics[bb= 10 130   1200   1300,height=7.5cm,width=8cm]{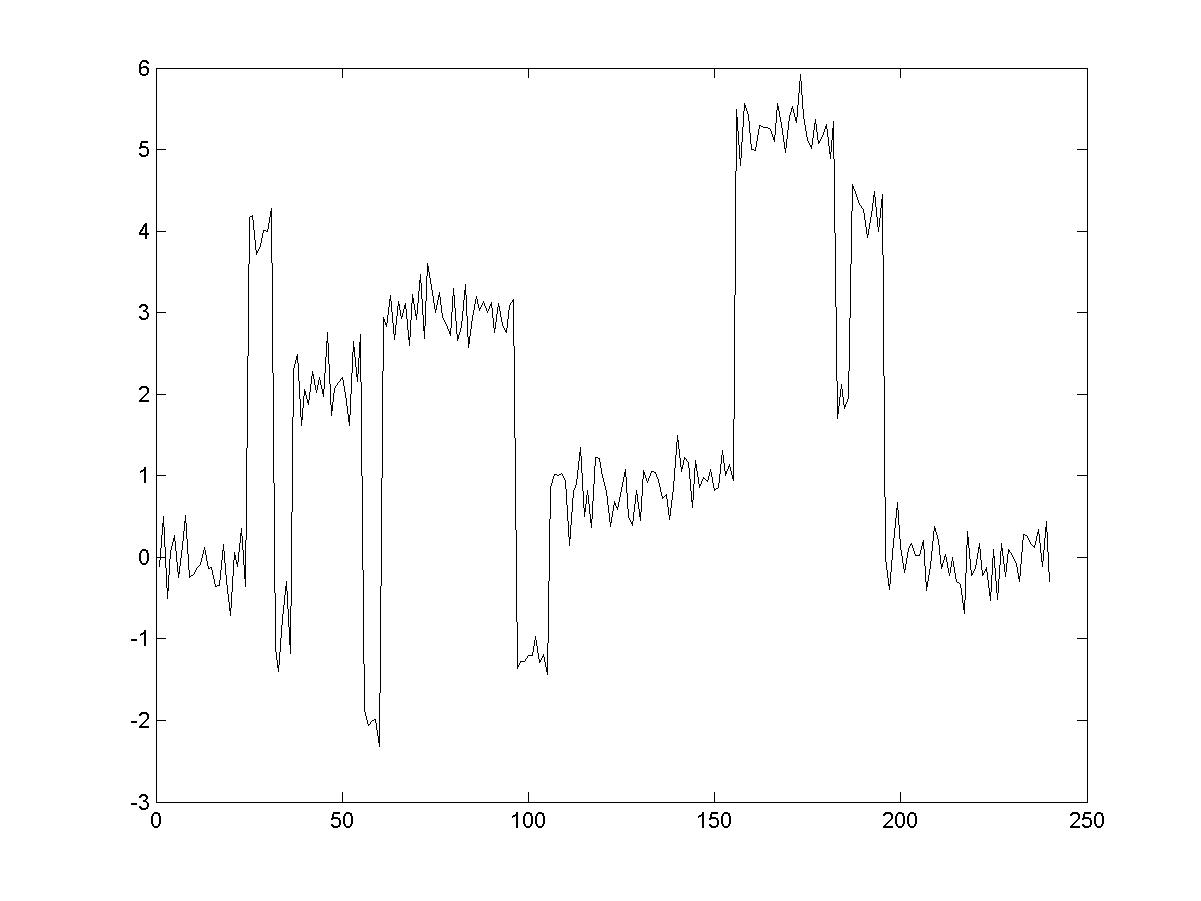}}
 \vspace{.2cm}
\centerline{
\includegraphics[bb= 10 130   1200   1300,height=7.5cm,width=8cm]{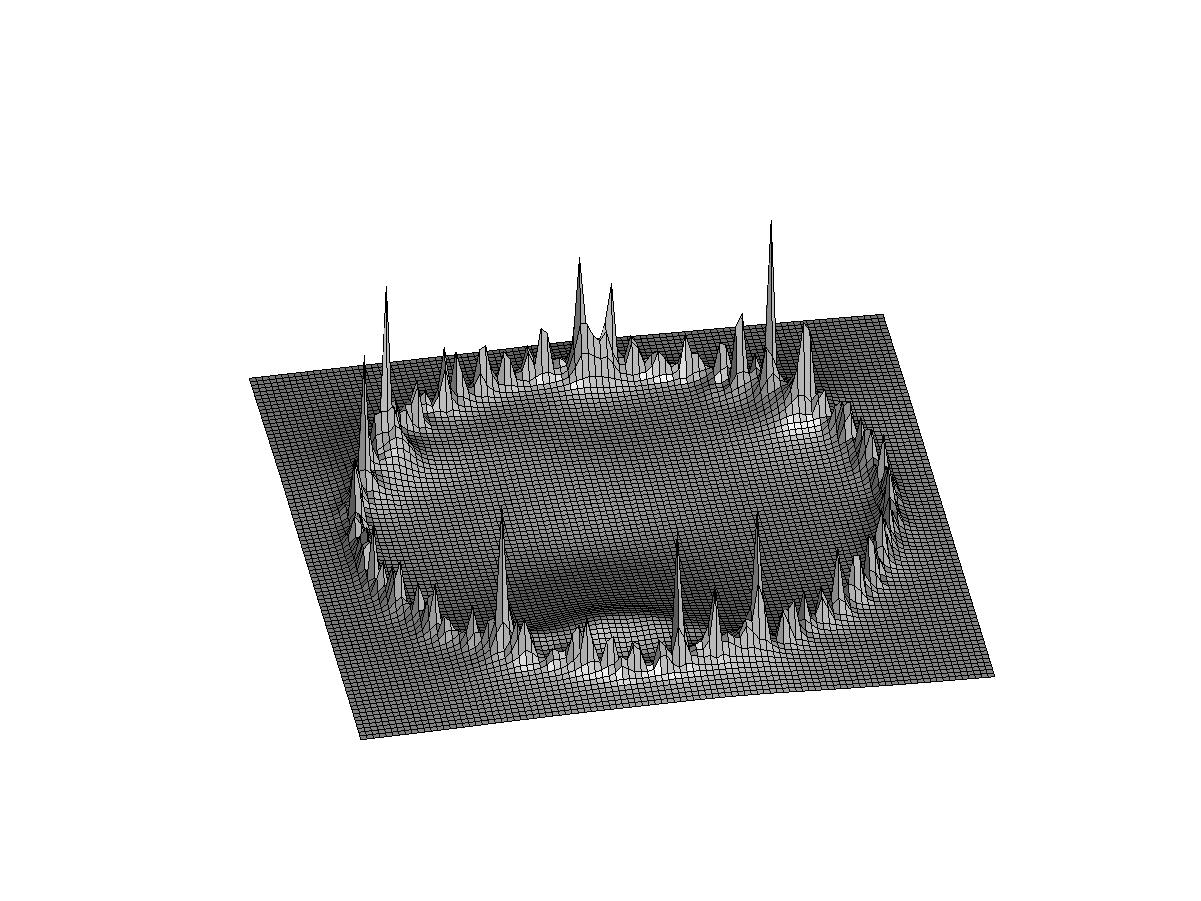}
\hspace{.5cm}
\includegraphics[bb= 10 130   1200   1300,height=7.5cm,width=8cm]{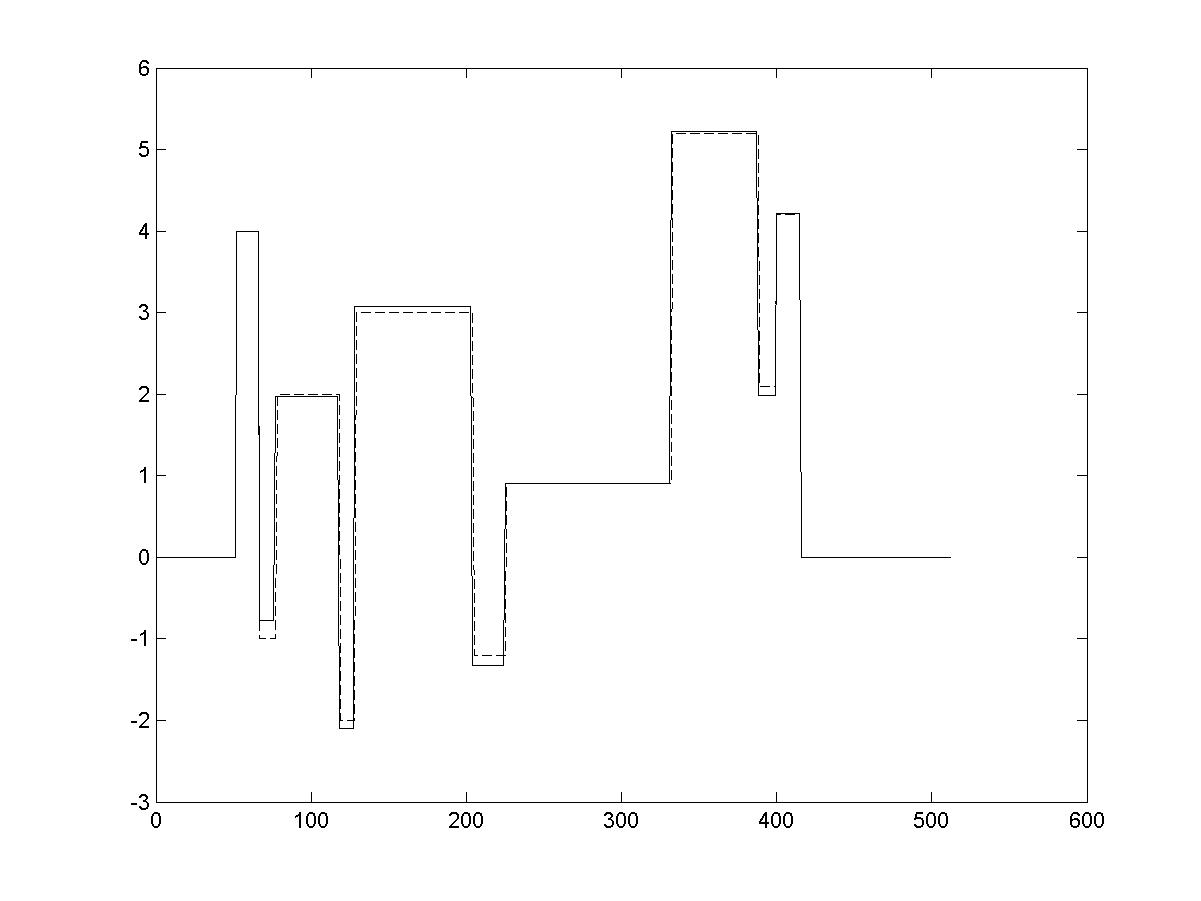}}
\end{center}
\caption{Top left: original function; top right: rough
estimate of $F(t)$ when the moments are affected by a
Gaussian noise with $SNR=7$. Bottom left: estimate of the
condensed density by the closed form approximation; bottom
right: reconstruction of the original function.}
\label{fig2}
\end{figure}

\begin{figure}[h]
\begin{center}
\hspace{1.7cm}\centerline{\includegraphics[bb= 10 130  1200   1300,height=7.5cm,width=8cm]{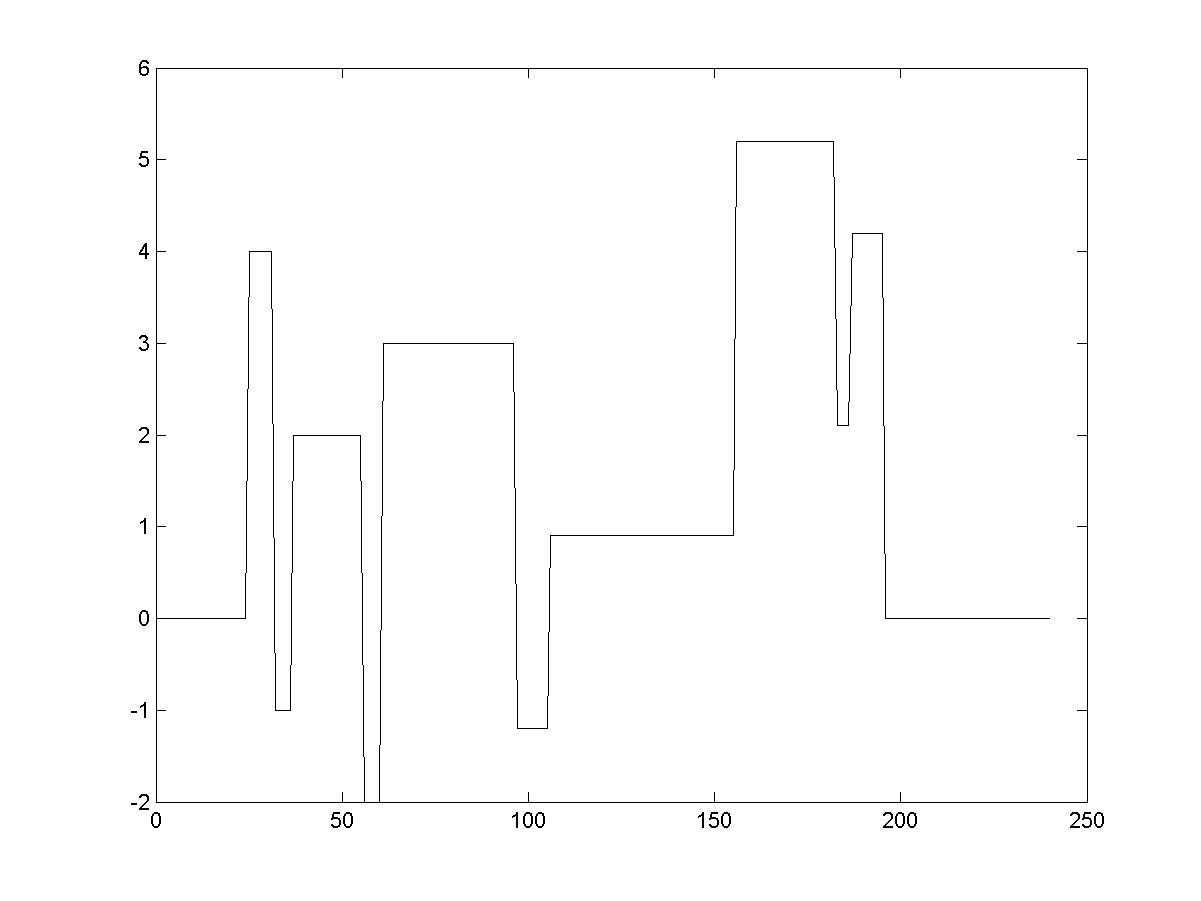}
\hspace{.5cm}\includegraphics[bb= 10 130   1200   1300,height=7.5cm,width=8cm]{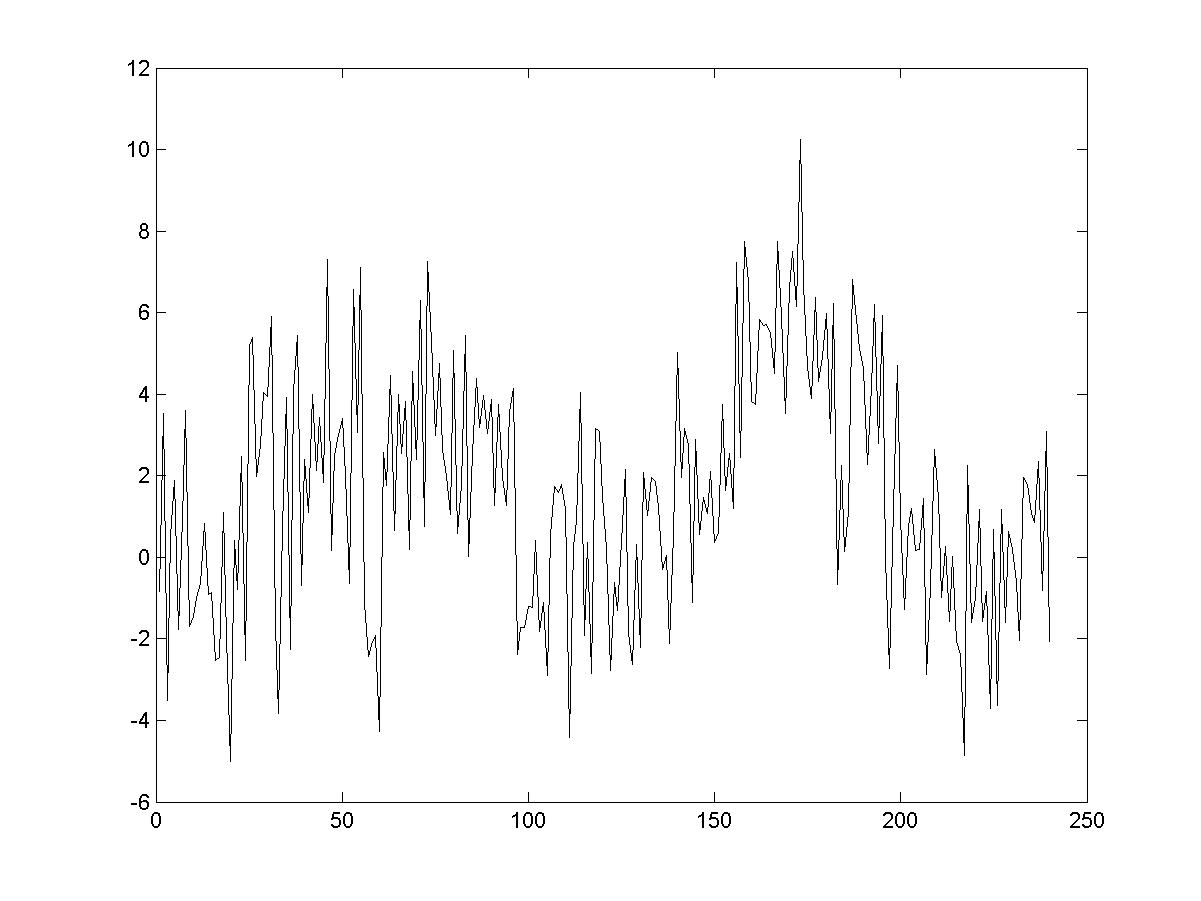}} \vspace{.2cm}
\hspace{1.7cm}\centerline{\includegraphics[bb= 10 130   1200   1300,height=7.5cm,width=8cm]{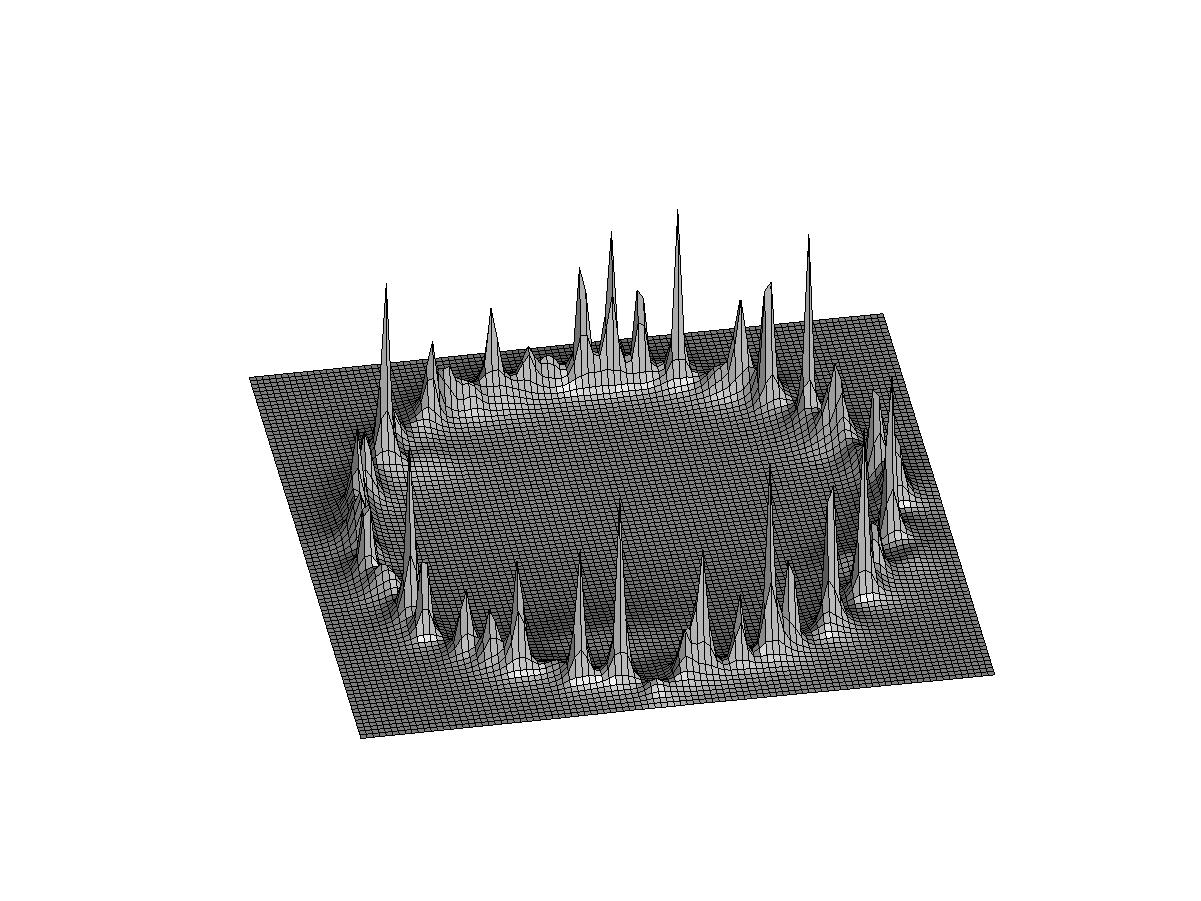}
\hspace{.5cm}\includegraphics[bb= 10 130   1200   1300,height=7.5cm,width=8cm]{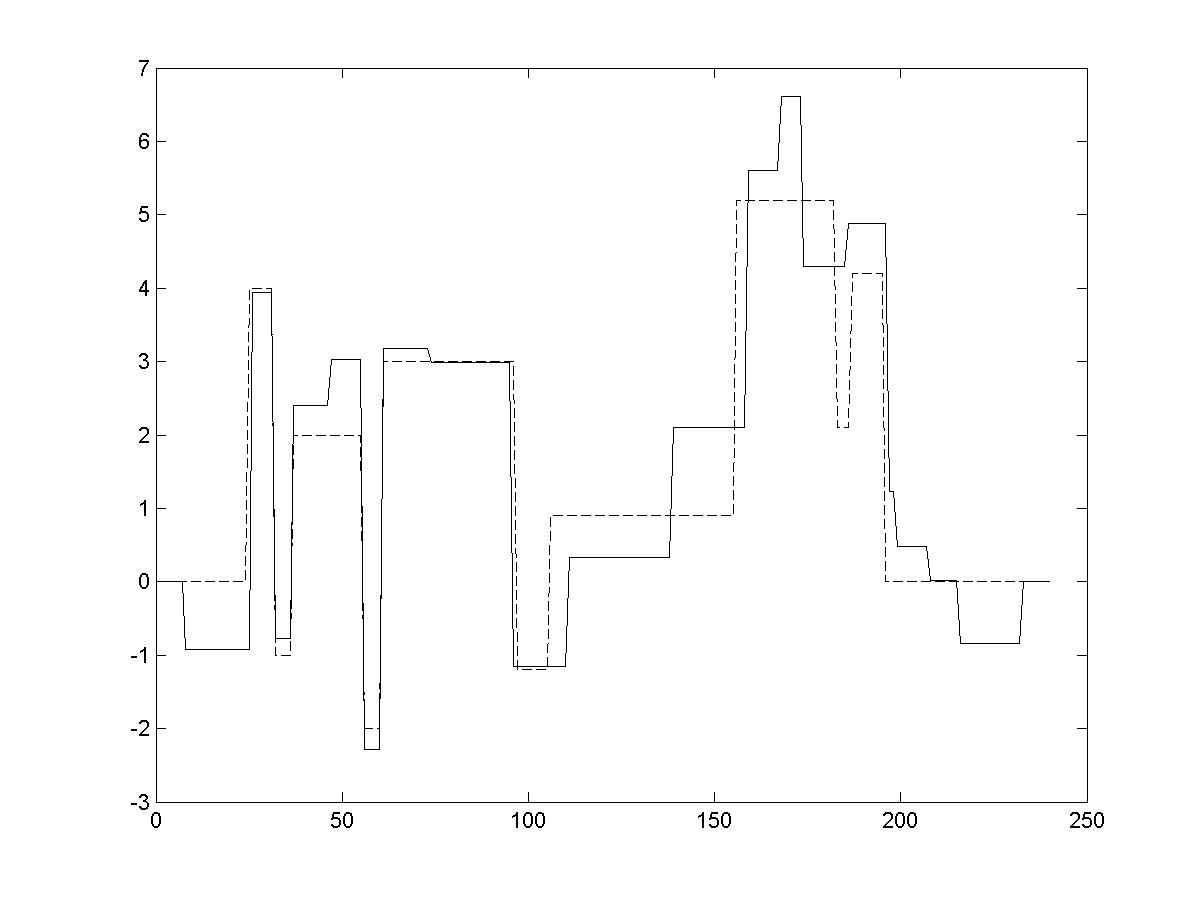}}
\end{center}
\caption{Top left: original function; top right: rough
estimate of $F(t)$ when the moments are affected by a
Gaussian noise with $SNR=1$. Bottom left: estimate of the
condensed density by the closed form approximation; bottom
right: reconstruction of the original function.}
\label{fig3}
\end{figure}

The method was applied to an example considered in
\cite{maba1} where comparisons with other methods were also
reported. In the top left part of fig.\ref{fig2} the
original function $F(t)$ is plotted. In the top right the
rough estimate of $F(t)$ when $SNR=7$ is reported where the
$SNR$ is measured as the ratio of the standard deviations
of $\{s_k\}$ and $\{\epsilon_k\}$. In the bottom parts the
condensed density and the reconstructed function
$\hat{F}(t)$ are plotted. Looking at the condensed density
we notice that the model is strongly identifiable,
therefore the estimation procedure outlined above was
applied. In fig.\ref{fig3} the same quantities as above but
with $SNR=1$ are plotted. In this case the model is
identifiable but not strongly therefore the clustering step
does not work. The number of complex exponentials used to
get the reconstruction plotted in fig.\ref{fig3} is
$\hat{p}=20$ and was found by trial and errors.

We notice that when $SNR=7$ we get an almost perfect reconstruction,
better than that reported in \cite{maba1}. When $SNR=1$ the
reconstruction quality is worse as expected but still comparable
with the one reported in \cite{maba1}.

\begin{table}[tbh]
\begin{center}
\begin{tabular}{|c|c|c|c|c|}
\hline \hline
&$\p$&$bias(\hat{p})$&$s.d.(\hat{p})$&$MSE(\hat{p})$\\
\hline
$ $&5& 0.0000&    0.0000 &   0.0000\\
\hline\hline
&$\xi_j$&$bias(\hat{\xi}_j)$&$s.d.\hat{\xi}_j$&$MSE(\hat{\xi}_j)$\\
\hline
$j=1$& -0.2796 - 0.8606i& -0.0008 + 0.0001i&   0.0000& 0.0000\\
\hline
$j=2$&-0.1782 - 0.9344i& 0.0036 - 0.0010i &  0.0000& 0.0000\\
\hline
$j=3$&   0.3090 + 0.9510i &  0.0057 - 0.0064i &  0.0031& 0.0001 \\
\hline
$j=4$&   0.2487 + 0.9685i & -0.0058 + 0.0110 &  0.0019 &0.0002 \\
\hline
$j=5$& -0.4354 + 0.5993i & -0.0047 + 0.0054i &  0.0108& 0.0002 \\
\hline \hline
&$c_j$&$bias(\hat{c}_j)$&$s.d.(\hat{c}_j)$&$MSE(\hat{c}_j)$\\
\hline
$j=1$&6.0000 &  0.0440  &  0.1238 &   0.0173\\
\hline
$j=2$&3.0000 &  -0.0407 &   0.0688  &  0.0064\\
\hline
$j=3$&1.0000 &  0.0441  &  0.0736   & 0.0074\\
\hline
$j=4$&1.0000&  -0.6767  &  0.0808  &  0.4644\\
\hline
$j=5$& 20.0000& -0.1007 &   0.2574  &  0.0764\\
\hline
\end{tabular}
\end{center}
\caption{Statistics of the parameters $\hat{p}$,
$\hat{\xi}_j,j=1,\dots,\p$ and $\hat{c}_j,j=1,\dots,\p$ }
\label{tb2}
\end{table}

%prodotta da \distrF\chi2nc16.m

\end{document}